\documentclass{amsart}   
\usepackage{amssymb}
\usepackage[arrow,matrix,curve]{xy}
\newcommand\cinv{\operatorname{c-ind}}
\newcommand\falg{\cPwf{-\ZZ}\ZZ X}

\newcommand\cohL{{}^{\text{L}}\kern-.5pt H}
\newcommand\Psu[2]{\Psi_{\sus}^{#1}(#2)}
\newcommand\Ps[2]{\Psi^{#1}(#2)}
\newcommand\sus{\operatorname{sus}}
\newcommand\ialg{\mathcal A}
\newcommand\ridl{\mathcal I}
\newcommand\alg[1]{\let\tem#1%
	\ifx\tem1{{\mathcal A}_{\sigma}}\fi
	\ifx\tem2{{\mathcal I}_{\sigma}}\fi
	\ifx\tem3{{\mathcal A}_{\pa,\sigma}}\fi
	\ifx\tem4{{\mathcal I}_{\pa}}\fi
	\ifx\tem5{{\mathcal A}_{\pa}}\fi}
\newcommand\algD[1]{\let\tem#1%
	\ifx\tem f\cPwf{-\ZZ}\ZZ X\fi
	\ifx\tem0\ialg=\cPwf{-\ZZ}\ZZ X/\cPwf{\infty}{-\infty}X\fi
	\ifx\tem1\alg1=\cPwf{-\ZZ}{\ZZ}X/\cPwf{-\ZZ}{-\infty}X\fi
	\ifx\tem2\idl2=\cPwf{\infty}{\ZZ}X/\cPwf{\infty}{-\infty}X\fi
	\ifx\tem3\alg3=\cPwf{-\ZZ}{\ZZ}X/
		\bigl(\cPwf{\infty}\ZZ X+\cPwf{-\ZZ}{-\infty}X\bigr)\fi
	\ifx\tem4\idl4=\cPwf{-\ZZ}{-\infty}X/\cPwf{\infty}{-\infty}X\fi
	\ifx\tem5\alg5=\cPwf{-\ZZ}{\ZZ}X/\cPwf{\infty}{\ZZ}X\fi
	\ifx\tem r\ridl=\cPwf{\infty}{-\infty}X\fi}
\newcommand\idl[1]{{\alg{#1}}}
\newcommand\Hc{\operatorname{HC}}
\newcommand\Hd{\operatorname{HH}}
\newcommand\bPsi[2]{\Psi_b^{#1}(#2)}
\newcommand\bT{{}^b\kern-.5pt T}
\newcommand\cT{{}^c\kern-.5pt T^*}
\newcommand\ct{{}^c\kern-.5pt T}
\newcommand\bt{{}^b\kern-.5pt T}
\newcommand\cS{{}^c\kern-.5pt S^*}
\newcommand\cM{{}^c\kern-.5pt MX}
\newcommand\bcM{{}^c\kern-.5pt M_{\pa X}X}
\newcommand\cnu{{}^c\kern-.5pt\nu}
\newcommand\dcnu{{}^c\kern-.5pt\nu_{\pa}}
\newcommand\HoCh[2]{{\mathcal H}_{#1}(#2)}

\newcommand\bcoh{{}^b\kern-1pt H}
\newcommand\bcohL{{}^{\text{bL}}\kern-.5pt H}
\newcommand\abs{\operatorname{abs}}
\newcommand\rel{\operatorname{rel}}
\newcommand\inv{\operatorname{inv}}
\newcommand\Ind{\operatorname{Ind}}
\newcommand\Bind{\operatorname{Ind}_{\pa}}
\newcommand\tr{\operatorname{tr}}
\newcommand\str{\operatorname{str}}
\newcommand\Tr{\operatorname{Tr}}

\newcommand\rTr{\Tr_{\pa,\sigma}}
\newcommand\RTr{\Tr_{\text{R}}}
\newcommand\hTr{\widehat{\Tr}}
\newcommand\dTr{\Tr_{\sigma}}
\newcommand\hdTr{\widehat{\Tr}_{\sigma}}
\newcommand\iTr{\Tr_{\pa}}
\newcommand\hiTr{\widehat{\Tr}_{\pa}}
\newcommand\bTr{\operatorname{\overline{Tr}}}

\newcommand\tTr{\widetilde{\operatorname{Tr}}}
\newcommand\tD[1]{\widetilde{D}(#1)}

\newcommand\Lap{\Delta}
\newcommand\bH[1]{{\mathcal M}_1(#1)}
\newcommand\bF[2]{{\mathcal M}_{#1}(#2)}
\newcommand\Vc{{\mathcal V}_{\text{c}}}
\newcommand\Vb{{\mathcal V}_{\text{b}}}
\newcommand\In{\operatorname{In}}
\newcommand\ins[1]{{#1}^{\circ}}
\newcommand\Hom{\operatorname{Hom}}
\newcommand\igc[1]{\rho^{\infty}\Psi_{\text{c}}^{-\infty}(#1)}
\renewcommand\Sp[3]{{#2}_{\text{#1}}^{#3}}
\newcommand\SP{\operatorname{SP}}
\newcommand\fD[2]{\Phi^{#1}(#2)}

\newcommand\cff{\text{ff}_{\text{c}}}
\newcommand\pcs[4]{I^{#1}_{\text{phg}}(#2,#3;#4)}
\newcommand\vpcs[5]{I^{#1}_{#5,\text{phg}}(#2,#3;#4)}
\newcommand\cPf[2]{\Psi_{\text{c}}^{#1}(#2)}
\newcommand\cPwf[3]{x^{#1}\Psi_{\text{c}}^{#2}(#3)}
\newcommand\cP[3]{\Psi_{\text{c}}^{#1}(#2;#3)}
\newcommand\Diag{\varDelta}
\newcommand\cDiag{\varDelta_{\text{c}}}
\newcommand\chd{\Omega_{\text{c}}^{\ha}}

\newcommand\Kd{\text{Kd}}
\newcommand\mhd{\Omega^{-\ha}}
\newcommand\Bif{\operatorname{Bf}}
\newcommand\If{\operatorname{IF}}
\newcommand\etab{\overline{\eta}}
\newcommand\ASb{\overline{\operatorname{AS}}}
\newcommand\iF{\beta}
\newcommand\sF{\gamma}
\newcommand\tZ{\widetilde{Z}}

\newcommand\La[2]{#1[[#2]}
\newcommand\totimes{\otimes_{\text{top}}}
\newcommand\indlim{\operatornamewithlimits{ind~lim}}
\newcommand\Fo[1]{\Lambda^{#1}}
\newcommand\bFo[1]{{}^b\Lambda^{#1}}
\newcommand\cFo[1]{{}^c\Lambda^{#1}}
\newcommand\Tt{a_0 \otimes a_1 \otimes \ldots \otimes a_n}
\newcommand\SSs{\SS_{\sigma}}
\newcommand\SSb{\SS_{\pa}}

\newcommand\datver[1]{\def\datverp%
 {\par\boxed{\boxed{\text{Version: #1; Run: \today}}}}}
\datver{2.0A; Revised: 29-5-1996}
\newcommand\ha{\frac12}
\renewcommand\Re{\operatorname{Re}}

\newcommand\CC{\mathbb C}
\newcommand\NN{\mathbb N}
\newcommand\RR{\mathbb R}
\newcommand\ZZ{\mathbb Z}
\renewcommand\SS{\mathbb S}
\newcommand\pa{\partial}
\newcommand\supp{\operatorname{supp}}
\newcommand\ci{${\mathcal C}^{\infty}$}
\newcommand\dCI{\dot{\mathcal C}^{\infty}}
\newcommand\CI{{\mathcal C}^{\infty}}

\newcommand\CIc{{\mathcal C}^{\infty}_{\text{c}}}
\newcommand\CmI{{\mathcal C}^{-\infty}}
\newcommand\dCmI{\dot{\mathcal C}^{-\infty}}
\newcommand\Id{\operatorname{Id}}
\newcommand\Mand{\text{ and }}

\newcommand\Mas{\text{ as }}
\newcommand\Mat{\text{ at }}
\newcommand\Mfor{\text{ for }}
\newcommand\Mif{\text{ if }}
\newcommand\Min{\text{ in }}

\newcommand\Motwi{\text{ otherwise }}
\newcommand\Motwip{\text{ otherwise.}}
\newcommand\Mwith{\text{ with }}
\newcommand\Mwhere{\text{ where }}
\newcommand\ie{i.e.}

\newtheorem{theorem}{Theorem}
\newtheorem*{Itheorem}{Index Theorem}
\newtheorem*{SItheorem}{Simplified Index Theorem}
\newtheorem*{RItheorem}{Reduced Index Theorem}
\newtheorem{proposition}{Proposition}
\newtheorem{corollary}{Corollary}
\newtheorem{lemma}{Lemma}

\theoremstyle{definition}
\newtheorem{definition}{Definition}

\theoremstyle{remark}

\begin{document}

\title[Homology of pseudodifferential operators I]
{Homology of pseudodifferential operators I. Manifolds with boundary}

\author[Richard Melrose]{Richard Melrose}
\address{Department of Mathematics, MIT}
\email{rbm@math.mit.edu}
\author[Victor Nistor]{Victor Nistor}
\address{Department of Mathematics, Pennsylvania State University}
\email{nistor@math.psu.edu}

\thanks{The research of the first author was supported in part by the
National Science Foundation under grant DMS-9306389; the second author was
partially supported by NSF Young Investigator Award DMS-9457859 and a Sloan
research fellowship. This manuscript is available by ftp from {\bf
ftp-math-papers.mit.edu} and from {\bf http:{\scriptsize//}www-math.mit.edu%
{\scriptsize/}$\sim$rbm{\scriptsize/}rbm-home.html} or {\bf
http:{\scriptsize//}www.math.psu.edu{\scriptsize/}nistor{\scriptsize/}}.}

\dedicatory\datverp

\begin{abstract} The Hochschild and cyclic homology groups are computed for
the algebra of `cusp' pseudodifferential operators on any compact manifold
with boundary. The index functional for this algebra is interpreted as a
Hochschild $1$-cocycle and evaluated in terms of extensions of the trace
functionals on the two natural ideals, corresponding to the two filtrations
by interior order and vanishing degree at the boundary, together with the
exterior derivations of the algebra. This leads to an index formula which
is a pseudodifferential extension of that of Atiyah, Patodi and Singer for
Dirac operators; together with a symbolic term it involves the `eta'
invariant on the suspended algebra over the boundary previously introduced
by the first author.
\end{abstract}
\maketitle

\section*{Introduction\label{Sect.I}}
The Hochschild and cyclic homologies of the algebra of
pseudodifferential operators on a compact manifold (without boundary)
were described by Wodzicki, \cite{Wodzicki3,Wodzicki1,Wodzicki2};
Brylinski and Getzler \cite{Brylinski-Getzler1} clarified and
generalized the computation of the homology of the algebra of full
symbols, which is the quotient by the ideal of smoothing
operators. Wodzicki and then Guillemin
\cite{Guillemin1,Guillemin2,Wodzicki3,Wodzicki1} defined a `residue
trace functional', this generates the Hochschild cohomology in
dimension $0$ for both the full algebra and the algebra of complete
symbols. These computations are closely related to the Atiyah-Singer
index formula since the index can be expressed in terms of a
Hochschild pairing for the symbol algebra using the `index functional'
which is the image of the trace functional on the smoothing operators
under the boundary map in the long exact sequence in Hochschild
homology discussed by Wodzicki in \cite{Wodzicki2}.

For any compact manifold with boundary the algebra of `cusp'
pseudodifferential operators is a natural generalization of the usual
algebra of pseudodifferential operators in the boundaryless case. We
describe the Hochschild homology of this algebra and various of its ideals
and so deduce a pseudodifferential generalization of the
Atiyah-Patodi-Singer index theorem. Such a generalization was previously
discussed by Piazza (\cite{Piazza1}) but without explicit evaluation of the
component functionals, here interpreted as Hochschild cycles.

In case $X$ is a compact manifold without boundary, let $\Ps{\ZZ}X$ be the
algebra of $1$-step polyhomogeneous (\ie\ classical) pseudodifferential
operators of any integral order. The smoothing operators form an ideal,
$\ridl=\Ps{-\infty}X,$ and the quotient $\ialg=\Ps{\ZZ}X/\ridl$ is
identified, by choice of a quantization map, with the algebra of formal
sums of homogeneous functions of integral order on $T^*X\setminus0$ with
order bounded above and with a $\star$-product. A spectral sequence
argument, see \cite{Brylinski-Getzler1}, allows the (continuous) Hochschild
homology and cohomology to be identified with de Rham cohomology spaces
\begin{equation}
\Hd_k(\ialg)\simeq H^{2n-k}(S^*X\times\SSs),\ \Hd^k(\ialg)\simeq
H^k(S^*X\times\SSs)
\label{rbm.296}\end{equation}
where $S^*X=(T^*X\setminus0)/\RR^+$ is the cosphere bundle; the
additional circle factor can be viewed as a quotient (or compactification)
of the $\RR^+$ action. The pairing uses
the symplectic form. In particular $\Hd_0(\ialg)$ is one-dimensional; the
dual space $\Hd^0(\ialg)$ is spanned by the residue trace, $\RTr,$ of
Wodzicki. This can be defined as follows. If
$Q\in\Ps1X$ is self-adjoint, elliptic and positive its complex powers,
defined by Seeley \cite{Seeley1}, allow elements of $\Ps{\ZZ}X$ to be
regularized. In particular the function
\begin{equation}
Z(A;z)=\Tr(AQ^{-z}),\ A\in\Ps mX,\ m\in\ZZ,
\label{rbm.295}\end{equation}
is holomorphic in $\Re z>m+\dim X$ and extends to be meromorphic. The
residue trace $\RTr(A)$ is the residue at $z=0$ and is independent of the
choice of $Q;$ it only depends on the class of $A$ in $\ialg.$

Although $\log Q$ is not an element of $\Ps{\ZZ}X$ commutation with it
defines a derivation on the algebras, $D_{\log Q}:\ialg\ni A\longmapsto
[\log Q,A].$ This exterior derivation induces a map on Hochschild
cohomology which we denote $i_{\log Q};$ it corresponds under
\eqref{rbm.296} with exterior (\ie\ cup) product with a generator of
$H^1(\SSs).$

The regularized value, $\hTr(A),$ of $Z(A;z)$ at $z=0$ is a functional on
$\Ps{\ZZ}X$ which is not a trace functional. Fedosov's formula suggests that
the {\em index functional} be defined in terms of its boundary, as a
Hochschild cochain
\begin{equation}
\begin{gathered}
\If(A,B)=(b\Tr)(A,B)=\hTr([A,B])=(\pa\Tr)(A,B)\\
=\RTr(A[\log Q,B])=(i_{\log Q}\RTr)(A,B).
\end{gathered}
\label{rbm.298}\end{equation}
Here $\pa$ is the boundary map in the long exact sequence 
\begin{equation}
0\longrightarrow\Hd^0(\ialg)\longrightarrow\Hd^0(\Ps{\ZZ}X)
\longrightarrow\Hd^0(\ridl)\overset{\pa}\longrightarrow\Hd^1(\ialg)
\longrightarrow\cdots
\label{rbm.297}\end{equation}
which arises from the short exact sequence $0\longrightarrow
\ridl\longrightarrow \Ps{\ZZ}X\longrightarrow\ialg\longrightarrow 0$ since
$\ridl$ is H-unital (see \cite{Wodzicki4}). Thus the index functional is
also the image under the boundary map of the class $\Tr$ on $\ridl.$ The
ideal $\ridl$ is homologous to a matrix algebra so
$\Hd^k(\ialg)\simeq\Hd^k(\Ps{\ZZ}X)$ except in dimension $1,$ where the
index functional in \eqref{rbm.298} is an additional generator of
$\Hd^1(\ialg).$ The other equalities in \eqref{rbm.298} show that the
index functional can be expressed in terms of Wodzicki's residue trace and
the derivation defined by $\log Q$ on the symbol algebra.

The identification of $\If$ as the `index functional' may be justified as
follows. By Morita 
equivalence the matrix algebras over the algebras under discussion have the same
Hochschild homology. If $A=[A_{ij}]$ is an $N\times N$ matrix in $\Ps{\ZZ}X$
which is elliptic, and hence Fredholm as an operator on $\CI(X)^N,$ then
its image $[A]$ in the matrix algebra over $\ialg$ is invertible (and 
conversely). The index formula of Fedosov then reduces to
\begin{equation}
\Ind(A)=\If([A],[A]^{-1}).
\label{rbm.299}\end{equation}
Since $\If$ reduces to the fundamental class in $H^0(S^*X),$ under the
isomorphism in \eqref{rbm.296}, the
Atiyah-Singer index theorem becomes the identification of the image in
$H^{2n-1}(S^*X)$ of the Hochschild cycle $[A]\otimes[A]^{-1}$ for $\ialg.$
Even though one might only be interested in operators of order $1,$ or even
of order $0,$ it is very convenient to discuss the formula \eqref{rbm.299}
in terms of the Hochschild homology for the algebra of operators of all
integral orders, since this is finite dimensional whereas the homology of
the operators of order $0$ is not.

It is this general discussion which we extend to the case of a compact
manifold with boundary, so that \eqref{rbm.299} becomes a
pseudodifferential index formula in that context, extending the
Atiyah-Patodi-Singer index theorem for Dirac operators.

The basic algebra considered here consists of the `cusp' pseudodifferential
operators on a compact manifold with boundary. It appears, at least
implicitly, in the work of Livingston \cite{Livingston1}. The properties of
the algebra are discussed more systematically in \cite{Melrose45}, which is
largely based on joint work of Rafe Mazzeo and the first author. This
algebra, denoted $\cPf{\ZZ}X,$ is 
determined by a choice of defining function, or more precisely by a
trivialization of the normal bundle to the boundary of the compact manifold
with boundary 
$X.$ A change of defining function gives an isomorphic algebra. A general
process of microlocalization allows $\cPf{\ZZ}X$ to be considered as the
quantization of an algebra of vector fields on $X$ which, with the
convention for coordinate systems that at a boundary point the first
function $x$ defines the boundary and its trivialization locally, is
spanned by $x^2\pa_{x},$ and the tangential vector fields $\pa_y.$ Thus, at
least informally, elements of $\cPf{\ZZ}X$ can be regarded as `symbolic
functions' of these vector fields, with smooth coefficients. If
$x\in\CI(X)$ is a global boundary defining function fixing the
trivialization of the normal bundle then this Lie algebra of `cusp' vector
fields is 
\begin{equation}
\Vc(X)=\left\{V\in\CI(X;TX);Vx\in x^2\CI(X)\right\}.
\label{rbm.300}\end{equation}
We look at the cusp calculus instead of the more familiar b-calculus to
avoid certain problems with completeness. The relationship between the two
algebras (which means in particular that the index formula we obtain is valid for
elliptic b-pseudodifferential operators) is briefly described in
Appendix~\ref{A.CAbC}. The b-calculus $\bPsi0X$ and the `corresponding'
cusp calculus $\cPf0X$ have isomorphic norm closures in the bounded
operators on $L^2.$ In \cite{Melrose-Nistor1}, the K-theory of the b-calculus is
described. Since the K-theory is given in terms of the norm completion the
conclusions are equally valid for the cusp calculus. At a \ci\ level the
`cusp commutation' $[x^2\pa_x,x]=x^2$ (or rather $[x^{-1},x^2\pa_x]=1$)
in place of the `b-commutation' $[x\pa_x,x]=x$ makes quite a lot of
difference as regards the computation of both Hochschild and cyclic
homology. Nevertheless methods similar to those used here can be applied
directly to the b-calculus.

The algebra $\cPf{\ZZ}X$ reduces to $\Ps{\ZZ}X$ in the boundaryless case and
consists in the general case of continuous linear operators on $\CI(X).$
The subspace $\dCI(X)\subset\CI(X)$ of functions vanishing to all orders at
the boundary is preserved by this action. Furthermore conjugation by any
complex power $x^z$ of a boundary defining function preserves the
algebra. Just as we consider the algebra of all integral order
pseudodifferential operators it is very helpful to expand the algebra to 
\begin{equation*}
\cPwf{-\ZZ}{\ZZ}X=\bigcup\limits_{k\in\ZZ}\bigcup\limits_{m\in\ZZ}x^k\cPf mX.
\end{equation*}
Here, and below, the arbitrary integral powers of $x$ are denoted
$x^{-\ZZ}$ to emphasize that this filtration is in the opposite direction
to the other, order, filtration. There is an ideal of smoothing operators
quite analogous to the boundaryless case, namely 
\begin{equation*}
\ridl=\cPwf{\infty}{-\infty}X=\bigcap\limits_{k\in\ZZ}\bigcap\limits_{m\in\ZZ}
x^k\cPf mX.
\end{equation*}
This is again an H-unital algebra, homologous to a matrix algebra, and
writing the quotient as $\ialg$ therefore leads to a long exact sequence in
Hochschild cohomology, as in \eqref{rbm.297}. It is again the case that a
matrix of elements of $\cPwf{-\ZZ}{\ZZ}X$ which is `fully elliptic', in the
sense that its image in the matrix algebra on $\ialg$ is invertible, is
Fredholm as an operator on $\dCI(X)^N.$ Its index is given by the same
`Fedosov formula' as the first part of \eqref{rbm.298} and \eqref{rbm.299} 
\begin{equation}
\Ind(A)=\dim \ker (A) - \dim \ker (A^*)=
\If([A],[A]^{-1})=(\pa\Tr)([A],[A]^{-1}).
\label{rbm.301}\end{equation}
However the direct analogue of the second part of \eqref{rbm.298} does not
hold. To derive a formula for the new index functional, $\If,$ we need to
consider two other ideals.

The smoothing algebra $\ridl$ is the residual ideal for the joint
filtration by order and boundary power. Each of these filtrations gives
rise to a separate ideal; we consider the quotient ideals in $\ialg$
\begin{gather}
\idl4=\cPwf{-\ZZ}{-\infty}X/\ridl,\ \cPwf{-\ZZ}{-\infty}X=
\bigcup\limits_{k\in\ZZ}\bigcap\limits_{m\in\ZZ}x^k\cPf mX,\\
\idl2=\cPwf{\infty}{\ZZ}X/\ridl,\ \cPwf{\infty}{\ZZ}X=
\bigcap\limits_{k\in\ZZ}\bigcup\limits_{m\in\ZZ}x^k\cPf mX,\\
\alg1=\ialg/\idl4,\ \alg5=\ialg/\idl2\Mand\alg3=\ialg/(\idl2+\idl4).
\label{rbm.302}\end{gather}
These algebras are related by the following diagram of short exact sequences
\begin{equation*}
\xymatrix{0\ar[dr]&&0\ar[d]&0\ar[d]\\
&\idl2+\idl4\ar[dr]&\idl2\ar[d]&\idl2\ar[d]\\
0\ar[r]&\idl4\ar[r]&\ialg\ar[r]\ar[d]\ar[dr]&\alg1\ar[r]\ar[d]&0\\
0\ar[r]&\idl4\ar[r]&\alg5\ar[d]\ar[r]&\alg3\ar[d]\ar[r]\ar[dr]&0\\
&&0&0&0.}
\end{equation*}

The two filtrations are delineated by maps, respectively the symbol map and
the indicial map. The first is quite analogous to the symbol map in the
boundaryless case, and reduces to it over the interior. The choice of a cusp
structure defines a modified tangent bundle such that $\Vc(X)=\CI(X;\ct X).$
This bundle is isomorphic to the usual tangent bundle. If $\cT X$ is the
dual bundle then the symbol map gives a short exact sequence for each $m$
\begin{equation}
\begin{gathered}
0\longrightarrow \cPf{m-1}X\hookrightarrow \cPf mX
\overset{\sigma_m}\longrightarrow {\mathcal G}^m(\cT X)\longrightarrow 0\\
\Mwhere{\mathcal G}^m(\cT X)=\left\{a\in\CI(\cT X\setminus0);\text{
homogeneous of degree }m\right\}.
\end{gathered}
\label{rbm.303}\end{equation}
Similarly, if $\Psu m{\pa X}$ is the `suspended algebra' of pseudodifferential
operators as discussed in \cite{Melrose46}, consisting of the
pseudodifferential operators on $\RR\times\pa X$ which are
translation-invariant in $\RR$ and have convolution kernels vanishing
rapidly with all derivatives at infinity, then the indicial map gives a
short exact sequence for each $m$ 
\begin{equation}
0\longrightarrow x\cPf mX\hookrightarrow\cPf mX\overset{\In}\longrightarrow
\Psu m{\pa X}\longrightarrow 0.
\label{rbm.304}\end{equation}

The choice of a quantization and an operator $Q\in \cPf1X$ which is elliptic,
positive and selfadjoint gives isomorphisms
\begin{equation}
\idl2\simeq\dCI(\cS X)[q]],\ \alg1\simeq x^{-\ZZ}\CI(\cS X)[q]]
\label{rbm.306}\end{equation}
where $q=\sigma_1(Q)$ is the symbol of $Q$ and $\cS
X=(\cT X\setminus0)/\RR^+$ is the sphere bundle associated to $\cT X;$ here
$[q]]$ stands for Laurent series in $q^{-1},$ \ie\ for each element there
is an upper bound on the powers of $q$ but no lower bound.
The product is a $\star$-product induced by the
quantization. Similarly the choice of the boundary defining function $x$
and a normal fibration of $X$ at $\pa X$ induces isomorphisms
\begin{equation}
\idl4\simeq\Psu{-\infty}{\pa X}[[x],\ \alg5\simeq\Psu{\ZZ}{\pa X}[[x].
\label{rbm.307}\end{equation}
Here, the commutation between the Laurent series variable $x$ and the
suspended algebra of pseudodifferential (or smoothing) operators is through
$[x,\pa_t]=x^2.$ Thus the variable on the linear factor in $\RR_t\times\pa
X$ can be identified as $t=1/x.$ The double quotient is then doubly a Laurent
series algebra 
\begin{equation}
\alg3\simeq\CI(\cS_{\pa X}X)[[x,q^{-1}].
\label{rbm.308}\end{equation}

These identifications and filtrations allow us to compute the Hochschild
homology groups of all these algebras. The Hochschild cohomology groups can
be expressed in terms of appropriate cohomology groups and duality:
\begin{equation}
\begin{gathered}
\Hd^*(\idl2)\simeq H^*(\cS X\times\SSs),\\
\Hd^*(\alg1)                              
\simeq H^*_{\rel}(\cS X\times\SSs)\oplus
H^*(\cS_{\pa X}X\times\SSs),\\
\Hd^k(\idl4)\simeq \begin{cases}\CC&\Mfor k=0,1,\\ 0 &\Motwi\end{cases}\\
\Hd^k(\alg5)\simeq
\begin{cases}H^k(\cS_{\pa X}X\times\SSs\times\SSb)/\CC&\Mfor k=1,2,\\
H^k(\cS_{\pa X}X\times\SSs\times\SSb)&\Motwi\end{cases}\\
\Hd^*(\alg3)\simeq H^*(\cS_{\pa X}X\times\SSs\times\SSb)\Mand\\
\Hd^k(\ialg)\simeq \begin{cases}
\CC \oplus H^*_{\rel}(\cS X\times\SSs)\oplus H^*(\cS_{\pa X}X\times\SSs)/\CC 
&\Mfor k=1\\
H^*_{\rel}(\cS X\times\SSs)\oplus
H^*(\cS_{\pa X}X\times\SSs)
&\text{ otherwise.}\end{cases}
\end{gathered}
\label{hop1.326}\end{equation}
%

In particular we show that in each case
$\Hd^0$ is one-dimensional, with a generating `trace functional'. These
functionals are, for $\ridl$ the ordinary trace, for $\alg1,$ $\alg5$ and
$\alg3$ a `double residue trace' $\rTr$ which we define by analytic
continuation, for $\idl2$ the residue trace of Wodzicki denoted $\dTr$ and
for $\idl4$ a functional $\iTr$ induced from the trace functional $\bTr$
defined in \cite{Melrose46} on $\Psu{\ZZ}Y$ for any boundaryless manifold   
$Y$ (its definition is recalled below). 
Analytic continuation arguments give extensions of the last two
functionals to $\hdTr$ on $\alg1$ and $\hiTr$ on $\alg5.$ These are not
trace functionals. Rather, commutation with the operators $\log Q$ and $\log
x$ defines derivations and the Hochschild boundaries of these functionals
are given by 
\begin{multline}
\begin{aligned}
(\pa\iTr)(A,B)&=\hiTr([A,B])=\rTr(A[\log Q,B])
=(i_{\log Q}\rTr)(A,B)\Mand\\
(\pa\dTr)(A,B)&=\hdTr([A,B])=-\rTr(A[\log x,B])
=-(i_{\log x}\rTr)(A,B)
\end{aligned}\\ \forall\ A,B\in\ialg.
\label{rbm.309}\end{multline}
The ideals $\idl2$ and $\idl4$ are H-unital and, in \eqref{rbm.309}, $\pa\dTr$
and $\pa\iTr$ represent the images under the respective boundary maps of the
functionals $\dTr$ and $\iTr$ in the long exact sequences 
\begin{gather}
0\longrightarrow\Hd^0(\alg1)\longrightarrow\Hd^0(\alg3)
\longrightarrow\Hd^0(\idl2)\overset{\pa}\longrightarrow\Hd^1(\alg1)
\longrightarrow\cdots\\
0\longrightarrow\Hd^0(\alg5)\longrightarrow\Hd^0(\alg3)
\longrightarrow\Hd^0(\idl4)\overset{\pa}\longrightarrow\Hd^1(\alg5)
\longrightarrow \cdots.
\label{rbm.310}\end{gather}
These formul\ae\ are analogues of the index formula \eqref{rbm.298}. Indeed
the sum of them gives a homotopy invariant of invertible elements
of $\alg3$ which we call the `boundary index' 
\begin{equation}
\Bind(A)=\Bif(A,A^{-1}),\ \Bif(A,B)=(-i_{\log x}\rTr+i_{\log Q}\rTr)(A,B).  
\label{rbm.312}\end{equation}
This represents an obstruction to the lifting of an invertible element
$A\in\alg3$ to an invertible element of $\ialg.$

The index functional itself, defined as the image in $\Hd^1(\ialg)$ of the
trace functional on $\Tr$ in \eqref{rbm.297} in the boundary case, can be
expressed in terms of these same operations 
\begin{multline}
\If(A,B)=\ha(\hiTr+\hdTr)(B\tD{A}+\tD{A}B),\\ A,B\in\ialg,\ \tD A=[\log
Q-\log x,A].
\label{rbm.311}\end{multline}

This leads to a pseudodifferential generalization of the
Atiyah-Patodi-Singer index theorem in which the index is expressed as the sum
of four terms. If $A=[A_{ij}]$ is a matrix in $\cPf mX$ which is elliptic and
also translation-invariant near the boundary, just as is usually assumed in
the case of Dirac operators, then one of these terms can be arranged to
vanish. If it is further assumed that the indicial family $A_0=\In(A)$ is
normal (for Dirac operators it is selfadjoint) then a second term can be
arranged to vanish and the resulting pseudodifferential index formula takes
the form of the Atiyah-Patodi-Singer formula
\begin{gather}
\Ind(A)=\ASb(A)-\ha\eta(A_0),\\
\ASb(A)=\ha\RTr(A^{-1}[\log Q,A]+[\log Q,A]A^{-1}),\\
\eta(A_0)=2\bTr(A_0^{-1}[\log x,A_0])
\label{rbm.313}\end{gather}
where $Q=(A^*A)^{\ha},$ the term $\ASb(A)$ is purely symbolic and
$\eta(A_0)$ is the `eta' functional defined on $\Psu{\ZZ}{\pa X}$ in
\cite{Melrose46}.

In Section~\ref{Sect.FADC} the Hochschild homology of algebras of functions
on a compact manifold with boundary is discussed. This is used in
Sections~\ref{Sect.HH}, \ref{S.BS} and \ref{S.IQ} to analyze the Hochschild
homology of the algebras $\idl2,$ $\idl4,$ $\alg1$ and $\alg3.$ The various
trace functionals are introduced explicitly in Section~\ref{Sect.RF} by
analytic continuation and the exterior derivations on the algebras are
considered in Section~\ref{Sect.ED}; their relation to the boundary maps is
discussed in Section~\ref{Sect.EC}. This leads to a description of the
Hochschild homology of $\alg5$ and $\ialg$ in
Section~\ref{Sect.THOA}. Morita invariance is then discussed and formula
\eqref{rbm.312} for the boundary index is obtained in
Section~\ref{Sect.BI}. The representation \eqref{rbm.311} of the index
functional is derived in Section~\ref{Sect.I1C} and the index formula,
\eqref{rbm.313}, is discussed in Section~\ref{Sect.IF}. The cyclic homology
of the main algebras is deduced in Section~\ref{Sect.CH}, the extension to
the $\ZZ_2$-graded case is described in Section~\ref{S.BaG} and the
homology of some other related algebras is examined in
Section~\ref{Sect.OA}. The appendices are devoted to a description of the
relevant properties of the algebra of cusp pseudodifferential operators.

\section{Functions and de Rham cohomology\label{Sect.FADC}}

In this section we obtain some preliminary results including the homology of
various weighted de Rham complexes.

Consider a compact manifold with boundary $M$ and let $x\in\CI(M)$ be a
boundary defining function. We shall consider the algebra of smooth
functions on the interior of $M$ with Laurent series expansions at the boundary
\begin{equation*}
x^{-\ZZ}\CI(M)=\bigcup\limits_{j\in\ZZ}x^j\CI(M).   
\end{equation*}
Let $\dCI(M)\subset\CI(M)$ be the ideal of smooth functions vanishing to
infinite order at the boundary. Then, with ${\CI(\pa M)}[[x]=
\{\sum_{k=-n}^{\infty}f_k x^{k}\}$ denoting the
Laurent series decomposed in any local product decomposition, Taylor's
theorem gives a short exact sequence
\begin{equation}
0\longrightarrow \dCI(M)\longrightarrow x^{-\ZZ}\CI(M)\longrightarrow
\La{\CI(\pa M)}x\longrightarrow 0.
\label{rbm.187}\end{equation}

Recall, \cite{Wodzicki4}, that an algebra $\mathcal{I}$ is H-unital if 
its  $b'$ complex is acyclic; the differential $b'$ is defined in
\eqref{hop1.331} below. The following result is a refinement of the
same result for $\CI(M),$ which is well known.

\begin{lemma}\label{rbm.188} The algebra $\dCI(M)$ is H-unital.
\end{lemma}

\begin{proof} For any two manifolds $M_1$ and $M_2$ the topological tensor
product is 
\begin{equation*}
\CI(M_1)\totimes\CI(M_2)=\CI(M_1\times M_2).
\end{equation*}
Since $\dCI(M)\subset\CI(M)$ is a closed subspace it follows that the space of
$k$-chains for Hochschild homology is 
\begin{equation*}
\dCI(M)\totimes\dCI(M)\totimes\cdots\totimes\dCI(M)=\dCI(M^{k+1}),
\end{equation*}
the space of smooth functions, on the manifold with corners $M^{k+1},$
vanishing to infinite order at all boundary faces. Let $z^{(j)}$ denote the
variable in the $j$th factor, for $j=0,1,\dots,k.$

Choose $\phi\in\CIc(\RR)$ with $\phi(0)=1.$ With $x^{(j)}$ denoting a fixed
defining function for the boundary of $M$ but in the $j$th factor consider 
\begin{equation}
Ef(z^{(0)},z^{(1)},\dots,z^{(k+1)})=
\phi\big(\frac{x^{(0)}-x^{(1)}}{x^{(0)}+x^{(1)}}\big)
f(z^{(1)},\dots,z^{(k+1)}),\ f\in\dCI(M^{k+1}). 
\label{rbm.189}\end{equation}
The rapid vanishing of $f$ at all boundaries, together with the fact that
the singular coefficient is supported near the diagonal as $x^{(0)}$ or
$x^{(1)}$ approaches $0,$ shows that 
\begin{equation*}
E:\dCI(M^{k+1})\longrightarrow \dCI(M^{k+2}).
\end{equation*}
By inspection,
\begin{equation*}
b'E+Eb'=\Id.
\end{equation*}
\end{proof}

Denote by 
$$
\chi(f_0 \otimes f_1 \otimes \ldots \otimes f_n)=
(n!)^{-1} f_0 df_1 \ldots df_n
$$
the Hochschild-Kostant-Rosenberg map, where $f_0, \ldots , f_n \in A,$ for
a commutative algebra $A.$  The map $\chi$ satisfies $\chi \circ b=0,$
where $b$ is the Hochschild differential. For the algebra of smooth
functions on a compact manifold $\chi$ induces a morphism
$\Hd_n(\CI(M))\to\CI(M,\Lambda^n),$ see
\cite{Hochschild-Kostant-Rosenberg1,Loday-Quillen1}, where
$\CI(M;\Lambda^k)$ denotes the space of smooth $k$-forms on $M.$ For a
compact manifold with boundary the same multilinear differential operator
gives rise to maps
\begin{equation}
\begin{gathered}
\chi:\dCI(M^{k+1})\longrightarrow \dCI(M;\Lambda^k),\\
\chi:x_1^{-\ZZ}\cdots x_{k+1}^{-\ZZ}\CI(M^{k+1})\longrightarrow
x^{-\ZZ}\CI(M;\Lambda^k)\Mand\\
\chi:\CI((\pa M)^{k+1})[[x_1,\dots,x_{k+1}]\longrightarrow
\CI(M;\Lambda^k)[[x]\oplus\CI(M;\Lambda^{k-1})[[x]\frac{dx}{x^2}
\end{gathered}
\label{hop1.352}\end{equation}
where $x_i$ is a fixed defining function for the boundary of $M$ but on the
$i$th factor.

Also recall from \cite{Wodzicki4} that, whenever the kernel is H-unital a
short exact sequence of algebras induces a long exact sequence in
Hochschild homology. From this it follows directly that the maps in
\eqref{hop1.352} all induce isomorphism in Hochschild homology.

\begin{lemma}\label{rbm.190} The short exact sequence \eqref{rbm.187}
induces a long exact sequence in Hochschild homology which decomposes into
the short exact sequences
\begin{multline}
0\longrightarrow\dCI(M;\Lambda^k)\longrightarrow 
x^{-\ZZ}\CI(M;\Lambda^k)\longrightarrow\\
\CI(\pa M;\Lambda^k)[[x]\oplus\CI(\pa M;\Lambda^{k-1})[[x]\frac{dx}
{x^2}\longrightarrow 0.
\label{rbm.191}
\end{multline}
Consequently the Hochschild-Kostant-Rosenberg maps in \eqref{hop1.352}
give isomorphisms of the Hochschild homologies of the three algebras in
\eqref{rbm.187}.
\end{lemma}

Notice that we could just as well use the basis element $dx/x$ (or indeed
$dx)$ in the third space in \eqref{rbm.191}. We use $dx/x^2$ because we
think in terms of c- (cusp) geometry rather than b-geometry.

\begin{proposition}\label{rbm.192} The de Rham differential commutes with
the maps in \eqref{rbm.191} so induces a long exact sequence of de Rham
cohomology spaces 
\begin{multline}
\cdots\longrightarrow H^k_{\rel}(M)\longrightarrow 
H^k_{\abs}(M)\oplus H^{k-1}(\pa M)\longrightarrow\\
H^{k}(\pa M)\oplus H^{k-1}(\pa M)\longrightarrow\cdots
\label{rbm.193}\end{multline}
which is just the direct sum of the Meyer-Vietoris sequence for $M$
relative to its boundary with the identity on $H^{k-1}(\pa M).$
\end{proposition}

\begin{proof}
The de Rham cohomology of the complex $\dCI(M;\Lambda^*)$ is isomorphic to
the cohomology of $M$ relative to its boundary by the de Rham theorem, since
it contains, and can be smoothly retracted onto, the complex of forms with
compact support in $M\setminus\pa M.$ Similarly the cohomology of the third
complex can be retracted onto the image of 
\begin{equation*}
\CI(\pa M;\Lambda^k)\oplus\CI(\pa M;\Lambda^{k-1})\ni(\alpha,\beta)
\longmapsto \alpha +\beta\wedge\frac {dx}x\in x^{-\ZZ}\CI(M;\Lambda^k),
\end{equation*}
where the Laurent series is decomposed using a local product decomposition
of $M$ near $\pa M.$ It follows that the cohomology is $H^{k}(\pa M)\oplus
H^{k-1}(\pa M).$ For the middle space a similar filtration argument shows
that the cohomology retracts onto the image of 
\begin{equation*}
\CI(M;\Lambda^k)\oplus\CI(\pa M;\Lambda^{k-1})\ni(\gamma,\beta )
\longmapsto \gamma +\beta \wedge\frac{dx}x\in z^{-\ZZ}\CI(M;\Lambda^k).
\end{equation*}
Again by the de Rham theorem the cohomology of the first summand is the
absolute cohomology of $M.$ Thus, the cohomology groups are as indicated in
\eqref{rbm.193} and the long exact sequence arises from  the short 
exact sequence of complexes
\begin{equation}
0 \rightarrow\dCI(M;\Lambda ^k)\rightarrow 
\CI(M;\Lambda ^k)\oplus\CI(\pa M;\Lambda^{k-1})\rightarrow 
\CI(\pa M;\Lambda ^k \oplus \Lambda^{k-1}) \rightarrow 0,
\label{rbm.208}\end{equation}
where the second map is the sum of pull-back to the boundary and the
identity map on $\CI(\pa M;\Lambda^{k-1}).$ Since, again by a form of the
de Rham theorem on a compact manifold with boundary, \eqref{rbm.208}
induces the Meyer-Vietoris sequence, the lemma is proved.
\end{proof}

The second of the image spaces in \eqref{hop1.352} consists of the forms on
the manifold with boundary $M$ which are smooth in the interior and have
Laurent series expansions at the boundary. 
It is interesting to note that the homology of the de Rham complex for
these spaces
\begin{equation}
0 \longrightarrow x^{-\ZZ}\CI(M;\Lambda^{0})\overset d\longrightarrow\dots
\cdots\overset d\longrightarrow x^{-\ZZ}\CI(M;\Lambda^k)
\overset d\longrightarrow \dots\cdots
\label{rbm.209}\end{equation}
turns out to be the same as that of various natural subcomplexes defined in
terms of Lie algebras of vector fields on the manifold. In particular we
consider the tangent (b-) vector fields $\Vb(M)$ and the vector fields
corresponding to a choice of cusp structure $\Vc(M);$ these are described
briefly in Appendix~\ref{App.A} and Appendix~\ref{A.CAbC}. Both Lie algebras
define vector bundles, $\bt M$ and $\ct M,$ over $M$ which are each
isomorphic to the usual tangent bundle over the interior and are such that
for the corresponding form bundles
\begin{gather*}
\CI(M;\Lambda^k)\subset\CI(M;\bFo k)\subset\CI(M;\cFo k),\\
x^{-\ZZ}\CI(M;\Lambda^k)=x^{-\ZZ}\CI(M;\bFo k)=x^{-\ZZ}\CI(M;\cFo k).
\end{gather*}
Furthermore (essentially because of the origins in terms of Lie algebras of
vector fields) the de Rham differential restricts to give complexes
\begin{gather}
\cdots\overset d\longrightarrow\CI(M;\bFo k)
\overset d\longrightarrow 
\CI(M;\bFo{k+1})\overset d\longrightarrow\dots\Mand 
\label{hop1.353}\\
\cdots\overset d\longrightarrow\CI(M;\cFo k)
\overset d\longrightarrow 
\CI(M;\cFo{k+1})\overset d\longrightarrow\cdots.
\label{rbm.194}\end{gather}
We shall denote by $\bcoh^*(M)$ the `b-cohomology' of $M$ which is the
homology of the first of these complexes. The b-cohomology of a
compact manifold with boundary was computed by Nest and Tsygan
\cite{Nest-Tsygan1},
\begin{equation*}
\bcoh^k(M)\equiv H^k_{\abs}(M)\oplus H^{k-1}(\pa M).
\end{equation*}
This is a natural isomorphism at the level of cohomology
where the map onto the second factor is the evaluation of the coefficient
of $dx/x.$ Cohomologies of this type had been computed
in unpublished work of the first author and Rafe Mazzeo; these (easy)
computations show that the homology spaces of the c-de Rham
complex determined by any cusp structure on $M$ are also isomorphic to the
b-cohomology spaces.

Consider a real vector bundle, $V,$ over a compact manifold, $X,$ possibly
with boundary or corners. Removal of the zero section gives a
non-compact manifold, $V\setminus0_V$ which has a natural $\RR^+$-action
on it. We shall denote by $S^*V$ the quotient sphere bundle and consider
$S^*V\times\SSs\equiv(V\setminus 0_V)\big/\ZZ$
the quotient by the discrete group $\ZZ$ acting through multiplication by
$e^k,$ $k\in\ZZ,$ on the fibres. Here the identification is fixed by the
choice of a positive function $r\in\CI(V\setminus0_V)$ which is homogeneous
of degree $1.$

If $X$ has no boundary then the cohomology of the complex of formal sums of
homogeneous forms on $V\setminus0_V$ of integral degree bounded above,
\begin{equation}
\begin{gathered}
0\longrightarrow {\mathcal G}(V)\overset{d}\longrightarrow{\mathcal G}(V;\Fo1) 
\overset{d}\longrightarrow {\mathcal G}(V;\Fo2)\longrightarrow\cdots\\
\Mwhere{\mathcal G}(V;\Lambda^k)=\bigoplus_{j=-\infty}^*{\mathcal G}^j
(V;\Lambda^k)\Mand\\
{\mathcal G}^j(V;\Lambda^k)=\big\{u_j\in
\CI(V\setminus0_V;\Fo k)\text{ homogeneous of degree }j\big\}
\end{gathered}
\label{rbm.221}\end{equation}
is easily seen to be concentrated in homogeneity degree $0.$ That is, if
$\alpha_i$ is a basis of de Rham cohomology on $S^*X$ and $r$ is a positive
function as above then the cohomology of \eqref{rbm.221} is spanned by the
forms $\alpha_i$ and $\alpha_i\wedge\frac{dr}r.$ These forms map to a basis
of the cohomology of $S^*V\times\SSs.$ That is, we have proved the
following

\begin{lemma}\label{rbm.222} The cohomology of the de Rham complex
\eqref{rbm.221} is naturally identified, as a graded space, with
$H^*(S^*V\times\SSs).$
\end{lemma}

In case $X$ is a manifold with boundary consider the Laurent forms on this
vector bundle, that is the formal sums of smooth forms on the interior of
$V\setminus0_V$ which are homogeneous, with degree bounded above, and have
at most rational singularities at the boundary:
\begin{gather*}
0\longrightarrow x^{-\ZZ}{\mathcal G}(V)\overset{d}\longrightarrow
x^{-\ZZ}{\mathcal G}(V;\Fo1) 
\overset{d}\longrightarrow x^{-\ZZ}{\mathcal G}(V;\Fo2)\longrightarrow\cdots\\
x^{-\ZZ}{\mathcal G}(V;\Lambda^k)=\bigcup\limits_{k\in\ZZ}
x^{-k}\bigoplus_{j=-\infty}^*\big\{u_j\in
\CI(V\setminus0_V;\Fo k)\text{ homogeneous of degree }j\big\}.
\end{gather*}
Again the cohomology is seen to be concentrated in homogeneity
zero so, by a simple computation, the cohomology of the de Rham complex for
these spaces is naturally identified with
\begin{equation}
\bcoh^k(S^*V\times\SSs)
\label{rbm.216}\end{equation}
where $\SS^*V$ is the sphere bundle of $V,$ and $\SSs$ is the circle
arising from a homogeneous function on the fibres.

\section{Spectral sequence\label{Sect.HH}}

In this we compute the Hochschild homology groups for some of the algebras
described in the Introduction, namely those which can be readily deduced by
spectral sequence arguments similar to those of Brylinski and Getzler in
\cite{Brylinski-Getzler1}.

Recall that if $A$ is a topological algebra the Hochschild
homology groups $\Hd_*(A)$ are the homology groups of the 
complex $({\mathcal H}_n(A), b)$ where $\mathcal H_n(A)$
is the completion of $A^{\otimes n+1}$ and $b$ is the (complete) Hochschild
differential
\begin{equation}
\begin{gathered}
b'(\Tt)=\sum_{i=0}^{n-1} (-1)^ia_0\otimes\ldots\otimes 
a_i a_{i+1}\otimes\ldots\otimes a_n\\
b(\Tt)=b'(\Tt)
+(-1)^n a_na_0\otimes\ldots\otimes a_{n-1}. 
\end{gathered}\label{hop1.331}
\end{equation}

The first algebra we consider for a compact manifold with boundary is an
analogue of the symbol algebra in the boundaryless case.

\begin{proposition}\label{rbm.116} Set
$\alg1=x^{-\ZZ}\cPf{\ZZ}X/x^{-\ZZ}\cPf{-\infty}X,$ the algebra of `complete
symbols with arbitrary polynomial blow up at the boundary'. Then
\begin{equation*}
\Hd_*(\alg1)\simeq \bcoh^{2n-*}(\cS X\times\SSs).
\end{equation*}
\end{proposition}

\noindent
Here we have used the notation of \eqref{rbm.216} for the `homogeneous'
cohomology of a vector bundle over a manifold with boundary. Thus, if
$\cS X$ is the (projective) sphere bundle of the cusp cotangent bundle,
which is diffeomorphic to the usual one and has boundary $\cS_{\pa X}X,$ then 
\begin{multline*}
\Hd_*(\alg1)\simeq\bcoh^{2n-*}(\cS X)\oplus\bcoh^{2n-1-*}(\cS X)\simeq\\
\big(H^{2n-*}_{\abs}(\cS X)\oplus H^{2n-1-*}(\cS_{\pa X}X)\big)
\oplus\big(H^{2n-1-*}_{\abs}(\cS X)\oplus H^{2n-2-*}(\cS_{\pa X}X)\big).
\end{multline*}

\begin{proof} This proceeds to a large extent as in the boundaryless case 
\cite{Brylinski-Getzler1}. Consider the complex $(\HoCh n{\alg1},b)$ which
computes the {\em continuous} 
Hochschild homology of $\alg1.$ The spaces ${\mathcal H}_n(\alg1)$
are obtained from $\alg1^{\otimes n+1}$ by topological completion. As
discussed in Appendix~\ref{App.A}, a choice of quantization map gives an
isomorphism 
\begin{equation}
\alg1=\indlim_{k,m\to\infty}
\left(x^{-m}\! \prod\limits_{j=-\infty}^{k}{\mathcal G}^j(\cT X)\right)
\label{hop1.328}\end{equation}
where ${\mathcal G}^j(V)$ is given by \eqref{rbm.221}.
Topological completion therefore gives the chain spaces a similar structure  
\begin{equation}
{\mathcal H}_{n-1}(\alg1)=
\indlim_{k,|\alpha|\to\infty}
\left(x^{-\alpha} \! \!\prod\limits_{j=-\infty}^{k} 
{\mathcal G}^j_n(\cT X) \right)
\label{hop1.327}\end{equation}
where $x_l$ is the boundary defining function in the $l$th factor,
$\alpha=(\alpha_1,\ldots,\alpha_n),$ $|\alpha|=\alpha_1+\cdots+\alpha_n,$
is a multiindex, $x^\alpha=x_1^{\alpha_1} x_2^{\alpha_2} \ldots
x_n^{\alpha_n},$ and ${\mathcal G}_n^j(V)\subset\CI((V\setminus0_V)^n)$ is
the subspace consisting of the 
functions which are finite sums of terms homogeneous in each factor with
total degree $j.$ The filtration of these spaces, given by `the
total order' is independent of the quantization used, so the corresponding
graded spaces are naturally defined and can be identified with the terms in
\eqref{hop1.327}.

The product on $\alg1$ is filtered with respect to the filtration by the
order, $k,$ in
\eqref{hop1.328}. It therefore gives rise to a spectral sequence as in
\cite{Brylinski-Getzler1}. The $E^1$ term of this 
spectral sequence computes the Hochschild homology of the
graded-commutative algebra $x^{-\ZZ}{\mathcal G}(\alg1)$ formed by the
$x^{-k}{\mathcal G}^j(\cT X).$  Corollary~\ref{rbm.190} shows that the
Hochschild-Kostant-Rosenberg map induces an isomorphism on
Hochschild homology
$$
\Hd_i({\mathcal G}(\alg1)) = x^{-\ZZ}{\mathcal G}(\cT X;\Lambda^i),
$$
where the image is the same formal sum of homogeneous $i$-forms (with
homogeneity bounded above) on $V\setminus0_V,$ for $V=\ct X.$

The differential $d_1,$ can be computed as in \cite{Brylinski-Getzler1}, 
Theorem 1 and Proposition 1. It is therefore the extension by continuity of
the Poisson differential on the interior; we denote this extension still by
$\delta.$ The naturality of the cusp exterior algebra, bearing as it does the same
relationship to cusp vector fields as the usual exterior algebra does to
all vector 
fields, shows that $\delta$ is the Poisson differential with respect to the
(degenerate) Poisson structure defined by the singular symplectic form on
$\cT X.$ In canonical local coordinates near the boundary this form is just
\begin{equation*}
\omega_{\text{c}}=\frac{dx}{x^2}\wedge
d\lambda+\sum\limits_{j=1}^{n-1}dy_j\wedge d\eta_j.
\label{hop1.329}\end{equation*}
Since it is homogeneous and the singularity at the boundary is of
finite order, symplectic duality as discussed by Brylinski
(\cite{Brylinski1}) gives isomorphisms 
\begin{gather*}
*_{\omega}:x^{-\ZZ}{\mathcal G}^j(\cT X;\Fo k)
\longrightarrow x^{-\ZZ}{\mathcal G}^{j+n-k}(\cT X;\Fo{2n-k}),\\
*(f dx^I d\xi^J) = \pm x^{2} f dx^{J^{\complement}} d\xi^{I^{\complement}},
\end{gather*}
transforming $\delta$ to the de Rham differential. Here $I^{\complement}$
denotes the complement of the index set $I.$ These identifications provide
an isomorphism of complexes
\begin{equation*}
(E^1_{p,q}, \delta)\simeq
\left(x^{-\ZZ}{\mathcal G}^{n-q}(\cT X;\Lambda^{2n-p-q}), d\right).
\end{equation*}

Thus the $E^2$ term of the spectral sequence vanishes, except for $q=n,$
the dimension of the manifold $X.$ From Lemma~\ref{rbm.222}, as modified in
the case of a manifold with boundary to \eqref{rbm.216}
\begin{multline*}
E^2_{p-n,n}=\bcoh^{2n-p}(\cS X\times\SSs)\\
=\bcoh^{2n-p}(\cS X)\oplus\bcoh^{2n-p-1}(\cS X) d\log r\\
=H^{2n-p}(\cS X\times\SSs) 
\oplus H^{2n-p-1}(\cS_{\pa X}X\times\SSs)x^{-1}dx
\end{multline*}
As in the boundaryless case it follows
that all further differentials are zero, so the
spectral sequence degenerates at $E_2.$

The convergence of the spectral sequence to the Hochschild homology also
follows directly, as in the boundaryless case.
\end{proof}

Next we consider the corresponding relative object which is an ideal in
$\alg1.$

\begin{proposition}\label{rbm.117} The algebra of `complete symbols
vanishing rapidly at the boundary',
$\idl2=x^{\infty}\cPf{\ZZ}X/x^{\infty}\cPf{-\infty}X,$ is H-unital and
$$
\Hd_*(\idl2)\simeq H_{\rel}^{2n-*}(\cS X\times\SSs)
\simeq H^{2n-*}_{\rel}(\cS X)\oplus H^{2n-1-*}_{\rel}(\cS X).
$$
\end{proposition}

\noindent
Here the `relative cohomology' groups, $H^{*}_{\rel}(\cS X)$ are the
de Rham cohomology groups computed relative to the boundary.

\begin{proof} The computation of the Hochschild homology proceeds very much
as in the previous proposition. One can also deduce it by retracting the
algebra to have compact support in the interior and then directly applying
the results of \cite{Brylinski-Getzler1}.

The H-unitality of $\idl2,$ \ie\ the fact that the corresponding $b'$
complex is acyclic, follows much as in Lemma~\ref{rbm.188}. Indeed, the
quantization can be chosen so that $\Id$ has full symbol expansion $1.$ The
cut-off factor in \eqref{rbm.189} then behaves as the identity with respect
to the $\star$ product.
\end{proof} 

\section{Boundary sequence\label{S.BS}}

Consider the short exact sequence 
\begin{equation}
0\longrightarrow\idl2\longrightarrow\alg1\longrightarrow\alg3\longrightarrow 0
\label{rbm.223}\end{equation}
where $\alg3$ is therefore the algebra of `Laurent series in $x$ of complete
symbols at the boundary'.

\begin{proposition}\label{rbm.118} For the algebra
$\alg3=x^{-\ZZ}\cPf{\ZZ}X/
\bigl(x^{\infty}\cPf{\ZZ}X+x^{-\ZZ}\cPf{-\infty}X\bigr)$ 
\begin{equation}
\Hd_*(\alg3)\simeq H^{2n-*}(\cS_{\pa X}X\times\SSs\times\SSb).
\label{hop1.330}\end{equation}
\end{proposition}

\begin{proof} 
Again the same approach can be used as in the proof of
Proposition~\ref{rbm.116}. This symbol algebra consists of formal power
series in the boundary variable $x$ as well as the homogeneous variable in
the fibres of $\cT X.$ This results in the two circles in \eqref{hop1.330},
with the de Rham cohomology being generated by $1, d \log x, d \log r$ and
$d \log x \wedge d \log r$ over $H^*(\cS_{\pa X}X).$
\end{proof}

\begin{proposition}\label{rbm.219} The long exact sequence in Hochschild
homology arising from \eqref{rbm.223} is the 
long exact sequence for b-
de Rham cohomology, of the manifold $Y_X=\cS X\times\SSs\times\SSb,$ with
respect to its boundary; \ie\ it is the usual 
cohomology long exact sequence
together with the identity on $H^*(\pa Y_X):$
\begin{equation*}
\xymatrix{&\ar[r]&
\Hd_k(\alg1)\ar[r]\ar@{=}[d]&\Hd_k(\alg3)\ar[r]\ar@{=}[d]
&\Hd_{k-1}(\idl2)\ar[r]\ar@{=}[d]&\\
&\ar[r]
&{\begin{matrix}H^{2n-k}(Y_X)\\ \oplus\\ H^{2n-k-1}(\pa
Y_X)\end{matrix}}\ar[r]^{i^*}_{\Id}
&{\begin{matrix}H^{2n-k}(\pa Y_X)\\ \oplus\\ H^{2n-k-1}(\pa Y_X)\end{matrix}}
\ar[r]^{\imath\pa}_{0}
&H^{2n-k+1}_{\rel}(Y_X)\ar[r]
&}
\end{equation*}
\end{proposition}

We have denoted by $i^*=i^*_{\pa Y_X}$ the restriction to the boundary.
Also note that the boundary map has an extra factor of $\imath =
\sqrt{-1}.$

\begin{proof} Since $\idl2$ is H-unital as a topological algebra it follows
from the result of Wodzicki, \cite{Wodzicki4}, that \eqref{rbm.223} induces
a long exact sequence in Hochschild homology 
$$
\cdots \longrightarrow \Hd_p(\alg1)
\longrightarrow \Hd_p(\alg3) \stackrel{\pa}{\longrightarrow} 
\Hd_{p-1}(\idl2) \longrightarrow \Hd_{p-1}(\alg1) \longrightarrow
\cdots.
$$
Except for the connecting morphism $\pa,$ the maps here come from algebra
morphisms. Since these morphisms preserve the order, they give
morphisms of the corresponding spectral sequences. This identifies  
the composition 
\begin{multline*}
H^{2n-p}_{\rel}(\cS X \times \SSs) \simeq 
\Hd_{p}(\idl2) \longrightarrow \\
\Hd_{p}(\alg1) \simeq H^{2n-p}(\cS X\times\SSs) 
\oplus H^{2n-1-p}(\cS_{\pa X}X\times\SSs)x^{-1}dx
\end{multline*}
as $\iota^* \oplus 0,$ where $\iota^*$ is the natural map from relative to 
absolute cohomology. Similarly, the composition
\begin{multline*}
H^{2n-p}(\cS X\times\SSs) 
\oplus H^{2n-1-p}(\cS_{\pa X}X\times\SSs)x^{-1}dx\simeq \Hd_p(\alg1)
\longrightarrow\\
\Hd_p(\alg3) \simeq H^{2n-p}(\cS_{\pa X}X\times\SSs)
\oplus H^{2n-1-p}(\cS_{\pa X}X\times\SSs) x^{-1}dx
\end{multline*}
preserves the direct sum decomposition, is the identity on the second
component and the restriction to the boundary on the first component.

It remains to determine the connecting morphism $\pa.$
The results of \cite{Wodzicki4} identify this map in the following
way. The inclusion of the complex ${\mathcal H}(\idl2)$ into
the complex ${\mathcal L},$ the kernel of the projection 
${\mathcal H}(\alg1) \to {\mathcal H}(\alg3),$ is shown to
induce an isomorphism of homology groups (in our case this can also be seen
from the associated spectral sequences). This shows that $\pa$ can be
computed as the connecting morphism of the homology
exact sequence associated to the short exact sequence 
\begin{equation}
\label{rbm.201}
0 \longrightarrow {\mathcal L} \longrightarrow 
{\mathcal H}(\alg1) \longrightarrow {\mathcal H}(\alg3) 
\longrightarrow 0.
\end{equation}
To carry out this computation assume first that $X=Y\times[0,1]$
where $Y$ is a manifold without boundary. Put
$Z_i=\cS_{Y\times \{i\}}X,$ $i=0,1.$ This gives, by restriction, an 
isomorphism $\alg3 \simeq \alg3^{(0)} \oplus \alg3^{(1)}.$
The group $\Hd_p(\alg3)$ decomposes then 
as a direct sum of two naturally isomorphic groups
$$
\Hd_p(\alg3) \simeq H^{2n-p}(\cS_{\pa X}X\times\SSs \times \SSb)
\simeq H^{2n-p}(Z_0\times\SSs \times \SSb) \oplus
H^{2n-p}(Z_1\times\SSs \times \SSb)
$$ 
and accordingly, so does the connecting morphism
$$
\pa = \pa_0 \oplus \pa_1: \Hd_p(\alg3) \to \Hd_{p-1}(\idl2) \simeq 
H^{2n-p+1}_{\rel}(\cS X\times\SSs).
$$

For reasons of symmetry it follows that $\partial_1 = -\partial_0,$
so it is enough to compute $\partial_0.$ Consider the subalgebra 
$(\alg3)_{\inv}$ consisting of those operators that commute with
$x^2(1-x)^2\pa_x.$  The algebra $\alg3$ is generated as a topological algebra 
by the subalgebra $(\alg3)_{\inv}$ and $x^{-\ZZ}(1-x)^{-\ZZ}\CI([0,1]).$
The algebra $\alg3^{(0)}$ is isomorphic to $x^{-\ZZ}(\alg3)_{\inv}.$ This defines
a splitting $s_0: \alg3^{(0)} \to \alg1.$
We remark that the existence of this splitting is a 
feature of the product and does not hold for arbitrary $X.$ 

Choose a smooth function $\varphi$ on $[0,1],$ which vanishes in a neighborhood
of $1,$ and is $1$ in a neighborhood of $0.$ and define 
\begin{gather*}
s: {\mathcal H}(\alg3^{(0)}) \to {\mathcal H}(\alg3) \\
s(f_0\otimes f_1 \otimes \ldots \otimes f_n)=
s_0(f_0) \varphi \otimes s_0(f_1) \otimes \ldots \otimes s_0(f_n)\Mand\\
\sigma : {\mathcal H}(\alg3^{(0)}) \to {\mathcal L}\\
\sigma (f_0\otimes f_1 \otimes \ldots \otimes f_n)=
s_0(f_0) [\varphi, s_0(f_1)] \otimes s_0(f_2) \otimes \ldots \otimes s_0(f_n)
\end{gather*}
where ${\mathcal L}$ is as in \eqref{rbm.201} above.

The relation $\sigma = [b,s],$ where
$b$ is the Hochschild boundary, follows from the fact that
$s_0$ is an algebra morphism and the definitions. An immediate
consequence is that $\sigma$ is a chain map with $[\sigma,b]=0$ in terms
of graded commutators. By definition, the natural morphism $\sigma_*,$
induced by $\sigma,$ from the Hochschild
homology groups of $\alg3^{(0)}$ to the homology of
${\mathcal L},$ coincides with the connecting morphism $\pa.$ The explicit 
formula
for $\sigma$ shows that it preserves the filtrations of both complexes
(shifting orders by $-1$). Since both spectral sequences degenerate at
$E^2,$ it is then enough to compute the morphism 
$E^2_{p,q}(\sigma):E^2_{p,q} (\alg3^{(0)}) \to E^2_{p-1,q}({\mathcal L}).$

Let $\eta$ with $d\eta =0,$ be a closed $(2n-p)$-form 
on $Z_0 \times\SSs \times \SSb.$ 
The computations in the proof of Proposition~
\ref{rbm.116} shows that we can represent the class of $\eta$ as
$$
[\eta] \in H^{2n-p}(Z_0 \times \SSs \times \SSb) 
\simeq E^2_{p-n,n} \simeq \Hd_{p}(\alg3^{(0)}),
$$ 
by a cycle $z \in F_{p-n}{\mathcal H}(\alg3^{(0)}),$ 
$z \in (\alg3^{(0)})^{\otimes p+1}.$ 

Denote by $\overline z$ the image of $z$ in the quotient 
${\mathcal G}^{p-n}(\alg3^{(0)}).$
We can assume that $\overline z$ is antisymmetric in all variables, and 
maps, via the Hochschild-Rosenberg-Kostant map $\chi,$ to
the form $*\eta,$ $\chi(\overline z) = *\eta.$ Let $H_{\varphi}$ be the 
Hamiltonian vector field associated to $\varphi,$ 
$H_{\varphi}(f)=\{\varphi,f\}.$
Here, as usual, we have denoted  the Poisson bracket by $\{\, , \,\};$ it
is determined by the singular symplectic structure on $\cT X.$ The antisymmetry
of $\overline z$ gives
$$
\chi(\sigma(\overline z))=\imath i_{H_\varphi}(\chi(\overline z)).
$$
Taking into account the symplectic duality $*$ and the relation 
$$
i_{H_f}(* \eta)= *(df \wedge \eta)
$$
we obtain the desired description of $\sigma_*$
$$
\pa([\eta])=\sigma_*([\eta]) = [ d\varphi \wedge \eta ] = [d (\varphi \eta)].
$$
where the brackets mean `the cohomology class' of that form. Since 
$\varphi \eta$ restricted to $Z_0$ is $\eta$
it follows that $[d (\varphi \eta)]$ is exactly the image of $[\eta]$ in
the relative group, via the connecting morphism in (de Rham) cohomology.

For arbitrary manifolds $X,$ choose an open  tubular neighborhood 
$V \simeq \pa X \times [0,1)$ of $\pa X$ in $X.$
The operators in $\alg1$ with support in $V$ identify with a subalgebra of
$\alg1(\pa X \times [0,1]).$ The computation follows then from the naturality
of the exact sequences in de Rham and Hochschild homology, and from the
naturality of the isomorphisms in Propositions \ref{rbm.116}, 
\ref{rbm.117} and \ref{rbm.118}.
\end{proof}

\section{Indicial quotient\label{S.IQ}}

Now consider the more interesting, indicial, quotient. That is, set
$$\idl4=x^{-\ZZ}\cPf{-\ZZ}X/x^{\infty}\cPf{-\ZZ}X.$$
This is the algebra of `indicial operators of order $-\infty$ with
Laurent series coefficients in $x$' as described in \eqref{rbm.307}.
That is, $\idl4$ is filtered by the powers of $x,$ $x$
being given a negative degree,
$$x^{-p}\cPf{-\ZZ}X/x^{-p+1}\cPf{-\ZZ}X\equiv
x^{-p}\Psu{-\infty}{\pa X}.$$
This isomorphism to the suspended algebra (which we also call the indicial
algebra) on the boundary is discussed in
the appendix. This filtration has completion properties similar to the 
the order filtration of the algebra of pseudodifferential
operators.

\begin{proposition}\label{rbm.119} The algebra $\idl4$ of Laurent series in
the indicial algebra is H-unital and
$$
\Hd_*(\idl4)\simeq H^{1-*}(\SSb).
$$
\end{proposition}

\begin{proof}
We use the filtration by powers of $x$ and the arguments similar to those in 
the preceding propositions.

The indicial algebra of order $-\infty$ is the
topological completion 
\begin{equation*}
\Psu{-\infty}Y\simeq  {\mathcal S}(\RR) \hat{\otimes} \Psi^{-\infty}(Y), 
\end{equation*}
where $\Psi^{-\infty}(Y)$ is the algebra of smoothing operators
on $Y,$  ${\mathcal S}(\RR)$ is the algebra of Schwartz functions on the line
with the convolution product, and $\hat{\otimes}$ is the projective
tensor product. 
 
The algebra $\Psi^{-\infty}(Y)$ 
is a projective module over its `double'
$$\Psi^{-\infty}(Y) 
\hat{\otimes}\Psi^{-\infty}(Y)^{\text{op}}$$
and hence is homologous to a matrix algebra, \ie\ its homology is
concentrated in dimension $0,$ where it is $\CC,$ given by the trace. This
shows that the structure of the graded algebra is
$$
\Psi^{-\infty}(\pa X) \hat{\otimes}{\mathcal S}(\RR) \hat{\otimes} 
\CC[x,x^{-1}].
$$
where all tensor products are projective tensor products and all three
algebras  are nuclear.
Using the topological K\"{u}nneth formula for Hochschild homology
\cite{Karoubi1} we see that the algebra of smoothing
operators does not change the result and hence the $E^1$-term
of the spectral sequence is 
\begin{equation}
E^1_{p,q}
 \simeq \Hd_{p+q}({\mathcal S}(\RR) 
\hat{\otimes}\CC[x]) = x^{-p}{\mathcal S}(\RR;\Fo{p+q-1}).
\end{equation}
The commutation relation $[x^{-1}, x^2 \pa_x]=1$ makes $\RR \times \SSb$
into a symplectic manifold, and the differential $d_1$ coincides with the
Poison differential as in the proof of Proposition~\ref{rbm.116}. 
The argument continues in the same way and show that the only
non vanishing $E^2$-terms are $E^2_{p,1}\simeq
H_{\text{comp}}^{2-p}(\RR\times\SSb),$ for $p=0,1.$ 

Since the spectral sequence clearly degenerates at $E^2$ and convergence
follows directly, we find
$$
\Hd_*(\idl4)\simeq H_{\text{comp}}^{2-*}(\RR\times\SSb)\simeq H^{1-*}(\SSb).
$$

The spectral sequence associated to the same filtration and Lemma 
\ref{rbm.188} shows that the $b'$-homology vanishes. This proves
H-unitality which can also be seen directly as in Lemma~\ref{rbm.188}.
\end{proof}

For the fully residual algebra $\ridl=x^{\ZZ}\cPf{-\infty}X$ one concludes
as in the boundaryless case that the Hochschild homology is just $\CC$ in
dimension $0.$

\section{Residue functionals\label{Sect.RF}}

In this section traces on various of the algebras are identified. As in
the boundaryless case, see \cite{Guillemin2,Wodzicki3}, we consider
functionals which arise as the residue of the analytic continuation of
`zeta-type' functions.  Fix a boundary defining function $x\in\CI(X)$ and a
positive, elliptic and invertible element $Q\in\cPf1X.$ For example, $Q$
can be taken to be $(\Lap+1)^{\ha}$ where $\Lap$ is the Laplacian of an
exact cusp metric, as discussed in Appendix~\ref{A.CM}.

\begin{lemma}\label{rbm.157} Suppose $A=A(z,\tau)\in\cPwf pmX$ is a
holomorphic function on $\Omega\subset\CC^2,$ where the connected
open set $\Omega$ contains a set of the form $[R,\infty)\times[R,\infty),$
then the function $\Tr(A(z,\tau)x^zQ^{-\tau})$ is defined and holomorphic
in the component of $\{\Re z>-p+1\}\cap\{\Re\tau>m+\dim X\}\cap\Omega$
meeting $[R',\infty)\times[R',\infty)$ for $R'$ large, and extends to a
meromorphic function on $\Omega$ with at most simple poles on the surfaces
$z=-p+1-j,$ $j\in\NN_0$ and $\tau=m+\dim X-k,$ $k\in\NN_0.$
\end{lemma}

\begin{proof} The general fact underlying this result is described in
Lemma~\ref{rbm.242} in Appendix~\ref{A.CM}. The present result follows
by applying this to the kernel of $A(z,\tau)$ as a conormal density on $\Sp
cX2$ with respect to the lifted diagonal $\cDiag.$ The function $x$ lifted
from the left factor is not quite a defining function for the boundary
hypersurface (defined by blow up) which the lifted diagonal meets, but is a
product $\rho \rho '$ where $\rho$ is such a function and $\rho '$ is a
product of boundary defining functions for boundary hypersurfaces at which
the kernel vanishes to infinite order. Thus 
the additional factor of $(\rho ')^z$ can be absorbed into the kernel to
give a conormal distribution which is entire in $z.$ Then
Lemma~\ref{rbm.242} applies, with the poles in $z$ shifted as a consequence
of the fact that the resulting density has an additional factor of $x^{-2}.$
\end{proof}

We are interested in the following
generalized (and `double') zeta function for $A\in\cPwf pmX$
\begin{equation}
Z(A;\tau,z)=Z_{x,Q}(A;\tau,z)=
\ha\big(\Tr(Ax^zQ^{-\tau})+\Tr(AQ^{-\tau}x^z\big)).
\label{rbm.156}\end{equation}
Lemma~\ref{rbm.157}, applied twice, shows that 
$\tau zZ_{x,Q}(A;\tau,z)$ is holomorphic in a neighborhood of
$0\in\CC^2.$ We shall examine the three functionals defined by 
\begin{equation}
\tau zZ(A;\tau,z)=\rTr(A)+\tau\hiTr(A) +z\hdTr(A)+\tau^2W+\tau zW'+z^2W''
\label{rbm.227}\end{equation}
where $W,$ $W'$ and $W''$ are holomorphic near $0.$

It also follows that the function  
\begin{equation*}
\tZ(A;z)=\tZ_{x,Q}(A;z)=Z_{x,Q}(A;z,z)=
\ha\big(\Tr(Ax^zQ^{-z})+\Tr(AQ^{-z}x^z)\big)
\end{equation*}
has at most a double pole at $z=0.$ The double residue of
$\tZ(A;z)$ at $z=0,$ (i.e. the coefficient of $z^{-2}$)
is just $\rTr(A).$ Clearly the residue of $\tZ(A;z)$ at $z=0$ is
$\hdTr(A)+\hiTr(A).$

Although we use the symmetrized functional in \eqref{rbm.156} the `left'
version 
\begin{equation*}
Z^L_{x,Q}=\Tr(Ax^zQ^{-\tau})
\end{equation*}
is only slightly different.

\begin{lemma}\label{hop1.332} For any $A\in\cPwf pmX,$ $
Z_{x,Q}(A;z,\tau)-Z^L_{x,Q}(A;z,\tau)$ is regular at $\tau=0$ and $z=0.$
\end{lemma}

\begin{proof} The difference is given by the analytic continuation of
$\Tr(AC(z,\tau)x^zQ^{\tau})$ where
$C(z,\tau)=\ha(\Id-Q^{-\tau}x^zQ^{\tau}x^{-z}).$ Since $C(z,\tau)$ is entire
of order $0$ and vanishes at $\tau=0$ and at $z=0$ the result follows.
\end{proof}

\begin{lemma}\label{rbm.158} The double residue $\rTr(A),$ of
$\tZ(A;z)$ at $z=0,$ is a trace functional on the algebra $\falg;$
it is independent of the choice of $Q$ or $x,$ vanishes on
$\idl2+\idl4$ and therefore defines a trace functional on $\alg3.$
\end{lemma}

\begin{proof} We can consider the double residue of
the non-symmetrized function  
\begin{equation}
\tZ^L_{x,Q}(A;z)=Z^L_{x,Q}(A;z,z).
\label{rbm.243}\end{equation}
If $x'$ is another defining function for the boundary then
the analytic family $A(z,\tau)=A-A(x'/x)^z$ also vanishes at $z=0$ so Lemma
\ref{rbm.157} now shows that the difference takes the form
$$
\tZ^L_{x,Q}(A;z)-\tZ^L_{x',Q}(A;z)=
\Tr(Ax^zQ^{-z}- A(x')^zQ^{-z})=z\Tr\bigl(B(z,\tau)x^zQ^{-z}\bigr)
$$
where $B(z,\tau)$ is holomorphic near $0$ and takes values in $\cPwf pmX.$
Thus the difference has at most a simple pole at $z=0$ and it follows that
the functional $\rTr(A)$ does not depend on the choice of the boundary
defining function. A similar argument shows independence of the choice of
$Q.$

To prove that $\rTr$ is a trace, consider $A$ and $B$ in the
`full' algebra $\falg.$ Expanding the commutator $[Q^{-\tau}A,Bx^z]$ and
using the trace property of $\Tr$ we deduce the following relation for large
real values of $z$ and $\tau$
\begin{equation}
\begin{gathered}
Z^L_{x,Q}([A,B];z,\tau)=\Tr([A,B]x^zQ^{-\tau})\\
=-\Tr\bigl([Q^{-\tau},B]Ax^z + BQ^{-\tau}[A,x^z] + B[Q^{-\tau},x^z]A\bigr)\\
=-\Tr\bigl((Q^{\tau}[Q^{-\tau},B]A + Q^{\tau}BQ^{-\tau}[Ax^{-z},x^z] + 
B[Q^{-\tau},x^z]AQ^{\tau}x^{-z})x^zQ^{-\tau} \bigr).
\end{gathered}\label{rbm.228}
\end{equation}
Since each of the families $Q^{\tau}[Q^{-\tau},B]A,$
$Q^{\tau}BQ^{-\tau}[Ax^{-z},x^z]$ and $B[Q^{-\tau},x^z]AQ^{\tau}x^{-z}$
is holomorphic as a family of classical pseudodifferential
operators  of fixed order, Lemma~\ref{rbm.157} shows that \eqref{rbm.228}
is valid in the domain of holomorphy. Due to
the commutators, each of these families vanishes at $z=\tau=0$ so
$\rTr([A,B])=0.$
\end{proof}

Before deriving an explicit formula for $\rTr(A)$ consider the functional
$\hdTr.$ Recall (see Appendix~\ref{App.A}) that if $2X$ is a compact manifold
without boundary obtained by doubling $X$ across its boundary then   
\begin{equation}
\cPwf\infty\ZZ X\subset\Psi^{\ZZ}(2X)
\label{rbm.244}\end{equation}
is an ideal, consisting precisely of those elements which have Schwartz
kernels supported in $X\times X.$ The smaller ideal, with kernels supported
in the interior, i.e. in $\ins{X}\times\ins{X},$ is dense. The residue trace
of Wodzicki for $2X$ is defined on the latter space; we shall
denote it $\RTr.$ We shall denote by $\dTr$ the restriction of $\hdTr$ to
$\idl2.$

\begin{lemma}\label{rbm.229} The restriction, $\dTr,$ of $\hdTr$  to
$\idl2,$ is a trace functional which is the extension by continuity of
Wodzicki's residue trace, $\RTr,$ for the double of $X.$
\end{lemma}

\begin{proof} If $A\in\idl2$ then $Ax^z$ is entire with values in $\idl2,$
as a family of fixed order. It
follows that $Z_{x,Q}(A;\tau,z)$ is entire in $z.$ Thus we only need
consider the simpler function 
\begin{equation}
Z_Q(A;\tau)=\Tr(AQ^{-\tau}),\ A\in\idl2.
\label{rbm.230}\end{equation}
>From a simplified version of Lemma~\ref{rbm.157} this is meromorphic in
$\tau$ with the pole at $\tau=0$ necessarily simple. The residue of this
function at $0$ is therefore just $\hdTr(A)=\dTr(A),$ by definition from
\eqref{rbm.227}. That is, 
\begin{equation}
\tau Z_Q(A;\tau)=\dTr(A)+\tau U(\tau),\Mfor A\in\idl2
\label{rbm.231}\end{equation}
with $U(\tau)$ holomorphic near $0.$ 

Suppose that $A\in\Psi^{\ZZ}(2X)$ is in the smaller ideal,
i.e.  has Schwartz kernel supported in $\ins{X}\times\ins{X}.$ Then
\eqref{rbm.231} is precisely the definition of the residue trace, except
that $Q$ should be positive and elliptic in $\Psi^1(2X).$ However, the
ellipticity of $Q$ in the cusp calculus implies that for any given $B$ with
kernel supported in $\ins{X}\times\ins{X}$ there exists
$Q'\in\Psi^1(2X)$ which is positive and elliptic and such that
$Q'B-QB\in\Psi^{-\infty}(2X).$ Since such a smoothing term does not
contribute to the residue trace, it follows that 
\begin{equation}
\dTr(A)=\RTr(A),\ A\in\Psi^{\ZZ}(2X),\ \supp(A)\subset\ins{X}\times\ins{X}.
\label{rbm.232}\end{equation}
Since $\RTr$ is itself a trace functional, the continuity of $\dTr$ in the
topology of $\idl2$ completes the proof. Alternatively the trace property
can be seen directly by a simple computation similar to that in
\eqref{rbm.228}. The independence of choice of $Q$ follows as before.
\end{proof}

This identification of $\dTr(A)$ and $\RTr(A)$ allows us to derive an
explicit formula for the former from the formula for the latter in
\cite{Guillemin2} or \cite{Wodzicki3}, see also \cite{Brylinski-Getzler1}.
Namely, any element $A\in\Psi^{\ZZ}(2X)$ can be identified with its
Schwartz kernel, a conormal right density with respect to the diagonal in
$(2X)^2.$ The collar neighborhood theorem allows a neighborhood of the
diagonal to be identified with a neighborhood of the zero section of $TX.$
By dropping a smooth term the kernel can be assumed to have support in this
neighborhood. Fourier transformation on the fibres converts it to a smooth
density on $T^*X.$ The polyhomogeneity of the original operator is just the
statement that this density has a complete asymptotic expansion 
\begin{equation*}
\bigl(\sum\limits_{j=-\infty}^{m}a_j\bigr)\omega^n,\ a_j\text{ homogeneous of degree }j,
\end{equation*}
(with the $n$th power of the symplectic form used to remove the density factor).
Even though the terms in this expansion depend on the choice of normal
fibration of $(2X)^2$ around the diagonal, the residue trace is always
given by
\begin{gather*}
\RTr(A)=(2\pi)^{-n}\int_{S^*(2X)} a_{-n}\nu
\end{gather*}
where $\nu$ is the $(2n-1)$-form obtained by contraction of the $n$th power
of the symplectic form on $T^*(2X)$ with the radial vector field,
$\zeta\cdot\pa_{\zeta}$ from the $\RR^+$ action on the fibres.

>From Lemma~\ref{rbm.229} and the continuity of $\dTr$ we deduce a
similar formula for $\dTr(A)$ when $A\in\idl2.$ This formula can be rewritten
in terms of the choice of a full symbol map as described in the
Appendix, and already used in the proof of Lemma~\ref{rbm.157}.

\begin{proposition}\label{rbm.233} If $A\in\idl2$ and
$\sum_{j=-\infty}^{m}a_j$ is the full symbol expansion of $A$ arising from
a choice of normal fibration to the lifted diagonal then 
\begin{equation}
\dTr(A)=(2\pi)^{-n}\int_{\cS X} a_{-n}\cnu
\label{rbm.234}\end{equation}
where $\cnu$ is the (singular) volume form on $\cS X$ defined by
contraction of the $n$th power of the symplectic form with the radial
vector field on the fibres of $\cT X.$ 
\end{proposition}

Notice that $\cnu$ is of the form $x^{-2}\nu'$ where $\nu'$ is a positive
smooth volume form on $\cS X.$ Thus we can see directly that $\dTr$ extends
by continuity from $\idl2$ to $\cPwf p\ZZ X,$ provided $p>1.$ In particular
we can consider the function  
\begin{equation}
\dTr(x^zA),\ A\in\falg
\label{rbm.235}\end{equation}
and observe that it is holomorphic for $\Re z$ large with a meromorphic
continuation to $z\in\CC.$ 

Now, consider a general element $A\in\falg,$ so $A\in\cPwf pmX$ for some
$p$ and $m\in\ZZ.$ As already noted this has a Laurent series expansion at
the boundary, with respect to any normal fibration and boundary defining
function $x$ 
\begin{equation}
A\sim\sum\limits_{j=-\infty}^{-p} x^{-j}A_j,\ A_j\in\Psu{m}{\pa X}.
\label{rbm.249}\end{equation}
In terms of this expansion we can obtain a formula for $\rTr(A).$ Observe
that the elements of $\Psu m{\pa X}$ are pseudodifferential operators on
$\RR\times\pa X,$ which are translation invariant in the first variable. A
normal fibration of $(\pa X)^2$ around the diagonal therefore leads to a
complete symbol expansion for each $A_j$ of the form $\sum\limits_k
a_{k,j}|\omega_{\pa}^{n-1}dtd\xi|$ where $\omega_\pa$ is the symplectic form
on $T^*\pa X$ and $\xi$ is the variable in the dual line (dual to $t$).

\begin{lemma}\label{rbm.226}
The functional $\rTr(A)$ for $A\in\falg$ is given explicitly by
\begin{equation}
\rTr(A)=(2\pi)^{-n}\int_{\cS_{\pa X}X} a_{-n,-1}\dcnu 
\label{rbm.195}\end{equation}
where $a_{k,l}$ is the term homogeneous of degree $k$ in the asymptotic
expansion of the symbol of $A_l$ with respect to a normal fibration of
$X^2$ around the diagonal and $\dcnu$ is the measure obtained by
contracting the form $\omega_{\pa}^{n-1}d\xi$ with the radial vector field
on $\cT_{\pa}X=T^*{\pa X}\times\RR_{\xi}.$
\end{lemma}

\begin{proof} Comparing the various functionals that we have so far defined
it follows directly that $\rTr(A)$ is the residue at $z=0$ of the
meromorphic function in \eqref{rbm.235}. Thus \eqref{rbm.195} follows
directly from \eqref{rbm.234}.
\end{proof}

Similarly we can deduce an explicit formula for $\hdTr.$ Namely  
it is just the Hadamard regularization of the usual formula for
Wodzicki's residue trace,
\begin{equation}
\hdTr(A)=\lim\limits_{\epsilon\downarrow0}
\left((2\pi)^{-n}\int_{\cS X\cap\{x>\epsilon\}} a_{-n}\cnu
+(\log\epsilon)\rTr(A)+\sum\limits_{l>0}\gamma_l\epsilon^{-l}\right),
\label{rbm.236}\end{equation}
where the coefficients are chosen so that the limit exists.

\begin{lemma}\label{rbm.246} The functional $\hdTr$ vanishes on $\idl4$ and
so defines a continuous functional on $\alg1,$ it is independent of the
choice of $Q.$ If $x'=ax,$ with $0<a\in\CI(X),$ are two boundary defining
functions the difference of the two functionals is given by 
\begin{equation}
\hdTr(A;x')-\hdTr(A;x)=\rTr(A\log a),\ A\in\falg.
\label{rbm.247}\end{equation}
\end{lemma}

\begin{proof} Independence of the choice of $Q$ follows from the fact that
$\RTr$ does not depend on $Q.$ The relation \eqref{rbm.247} between the
functionals obtained by regularization with respect to different boundary
defining functions follows from the fact that it is given by the regular
value at $z=0$ of 
\begin{equation*}
\dTr(x^z(a^z-1)A)=\dTr(zx^zB(z))
\end{equation*}
where $B(z)=z^{-1}(A^z-1)A$ is entire. Since the residue of 
$\dTr(x^zB(z))$
at $z=0$ is $\rTr(B(0))$ the result follows.
\end{proof}

Next we reverse the roles of $x$ and $Q$ and proceed to consider the
functional $\hiTr,$ defined on $\falg$ by \eqref{rbm.227}. Since
$Z(A;\tau,z)$ is entire if $A\in\ridl,$ the functional $\hiTr$ is
well-defined on $\alg5,$ we shall write its restriction to $\idl4$ as
$\iTr.$  If $A\in\cPf{-\infty}X$ then $Z(A;\tau,z)$ is entire in
$\tau,$ as a meromorphic function of $z.$ It follows then that $\iTr(A)$ is
the residue at $z=0$ of $\Tr(Ax^z).$ Recall from \cite{Melrose46} that on
the indicial algebra $\Psu{-\infty}{\pa X}\simeq {\mathcal S}(\RR) 
\widehat{\otimes} \Psi^{-\infty}(\pa X)$ there is a trace functional
given by integration of the trace of the indicial family  
\begin{equation}
\bTr(B)=(2\pi)^{-1}\int_{\RR}\Tr\widehat{B}(\xi)d\xi.
\label{rbm.237}\end{equation}

\begin{lemma}\label{rbm.238} The functional $\iTr$ is a trace functional on
$\idl4,$ defined independently of choice of $x$ or $Q,$ and is fixed in
terms of the trace \eqref{rbm.237} on the indicial algebra through 
\begin{equation*}
\iTr(A)=\bTr(A_{-1}),\ A\in\idl4,\ A\sim\sum\limits_{j=-\infty}^{-p}x^{-j}A_j,\
A_j\in\Psu{-\infty}{\pa X}.
\end{equation*}
\end{lemma}

\begin{proof} This follows directly from the formula for the trace as the
integral of the Schwartz kernel over the diagonal. Lifting to $\Sp cX2$ this
becomes 
$$
\Tr(Ax^z)=\int_{\cDiag} x^zA_{\cDiag}=\int_X x^z B\nu
$$
where $\nu=x^{-2}\nu',$ with $\nu'$ a smooth positive density on $X$ and
$B\in x^{-\ZZ}\CI(X)$ represents the kernel restricted to the
diagonal. The residue at $z=0$ therefore arises from the term $B_{-1}$ in
the expansion 
\begin{equation*}
\iTr(A)=\int_{\pa X} B_{-1}\nu_{\pa},\ \nu'=dx\nu_{\pa}\Mat\pa X,\
B\simeq \sum\limits_{k=-\infty}^{-p} x^{-k}B_k,
\end{equation*}
Since $B_{-1}\nu_{\pa}$ is the restriction of the kernel of $A _{-1}$ to
$\{0\}\times\Diag_{\pa X}$ it follows that the trace is given by
\eqref{rbm.237} with $B=A_{-1}.$ Clearly this is independent of the choice
of $Q,$ and independence of the choice of $x$ follows as above.
\end{proof}

As shown in \cite{Melrose46} the functional $\bTr$ extends from
$\Psu{-\infty}{\pa X}$ to $\Psu{\ZZ}{\pa X}$ still as a trace
functional. This extension is given by regularization of
\eqref{rbm.237}. Namely, if $B\in\Psu{m}{\pa X}$ then   
\begin{equation}
h_p(\xi)=\Tr\left((\frac{\pa}{\pa\xi})^p\widehat B(\xi)\right),\ p\ge m+\dim X
\label{rbm.248}\end{equation}
is defined and the asymptotic expansion of the $(p+1)$-fold integral
\begin{equation}
\int_{-T}^T \int_{0}^{\xi_p}\cdots\int_{0}^{\xi_1}
h_p(\xi)d\xi d\xi_1\cdots d\xi_p\sim\sum\limits_{j<m+\dim X} c_jT^j
+P(T)\log T
\label{rbm.239}\end{equation}
fixes the coefficient $\bTr(B)=c_0/2\pi$ independent of $p.$ Here $P$ is a
polynomial and the independence of $p$ follows from the fact that
increasing $p$ changes the inner $p$-fold integral in \eqref{rbm.239} by a
polynomial and hence changes the full integral by a polynomial without
constant term, leaving $c_0$ unchanged. 

Let $Q_0=(D_t^2+\Lap_{\pa X}+1)^{\ha}$ be a positive and elliptic element
of $\Psu1{\pa X},$ given by the Laplacian of a metric on $\pa X.$ Then
consider the function  
\begin{equation}
\bTr(BQ_0^{-\tau}).
\label{rbm.240}\end{equation}
This is holomorphic for $\Re\tau$ large and has a meromorphic continuation
to the complex plane. Indeed, if $B\in\Psu m{\pa X}$ then the poles can
occur only at $\tau=-\dim X+\NN_0.$  

\begin{proposition}\label{rbm.241} For $B\in\Psu m{\pa X},$ $m\in\ZZ,$ the
residue at $z=0$ of the function in \eqref{rbm.240} is the functional
$\tTr(B)=\bTr([t,B])$ introduced in \cite{Melrose46} and the regularized
value is $\bTr(B),$ i.e.
\begin{equation*}
\tau\bTr(BQ_0^{-\tau})=\tTr(B)+\tau\bTr(B)+\tau^2 U'(\tau)
\end{equation*}
with $U'(\tau)$ holomorphic near $\tau=0.$
\end{proposition}

\begin{proof} For $p$ large, let $h_p(\xi;\tau)$ be the function given by
\eqref{rbm.248} with $B$ replaced by $BQ_0^{-\tau}.$ This is entire in
$\tau$ and for $\tau$ near zero the asymptotic expansion as
$|\xi|\to\infty$ is also uniform with 
\begin{equation*}
h_p(\xi;\tau)\sim\sum\limits d_k^{\pm}(\tau)|\xi|^{m-p+n-1-\tau-k},\
\xi\to\pm\infty,
\end{equation*}
the coefficients being holomorphic, as is any finite order
remainder. Inserting this into \eqref{rbm.239} shows that $c_0(\tau)$ has
only a simple pole at $\tau=0$ with residue coming from $d_{k}^{\pm}$ for
$k=m+n$ and moreover the regular value at $\tau=0$ is just $\bTr(B).$ The
identification of the residue with $\tTr(B)$ follows from the results of
\cite{Melrose46}.
\end{proof}

One direct consequence of this formula is the identification of the residue
trace $\rTr$ in terms of the functional $\tTr$ on the indicial algebra.

\begin{corollary}\label{rbm.253} If $A\in\falg$ then the residue of
$\iTr(AQ^{-\tau})$ at $\tau=0$ is 
\begin{equation}
\rTr(A)=\tTr(A_{-1}).
\label{rbm.254}\end{equation}
\end{corollary}

\begin{corollary}\label{rbm.250} If the regularizer $Q$ is asymptotically
translation invariant with respect to a given fibration of $X$ near $\pa X,$
so that the expansion \eqref{rbm.249} for $Q$ becomes $Q\sim Q_0,$ then in
terms of the corresponding expansion for a general element $A\in\falg,$  
\begin{equation*}
\hiTr(A)=\bTr(A_{-1}).
\end{equation*}
\end{corollary}

\noindent
Indeed, $AQ\sim\sum_j x^{-j}A_jQ_0,$ so this follows from
Proposition~\ref{rbm.241}.

Note that $\hiTr$ is independent of the choice of $x$ but depends on the
choice of $Q.$ Using this corollary it is straightforward to deduce the
formula for a general regularizer $Q.$ If $Q'$ and $Q$ are both positive
and elliptic then the formal difference of logarithms
\begin{equation*}
\log(Q'/Q)=\frac{d}{d\tau}\left((Q')^{\tau}Q^{-\tau}\right)\big|_{\tau=0}
\end{equation*}
is a well-defined element of $\cPf0X$ which satisfies the cocycle condition
\begin{equation}
\log(Q''/Q')+\log(Q'/Q)=\log(Q''/Q).
\label{rbm.251}\end{equation}

\begin{lemma}\label{rbm.252} If $\hiTr(A;Q)$ denotes the value of the
functional $\hiTr$ on an element $A\in\falg$ given by a regularizer $Q$ then 
\begin{equation*}
\hiTr(A;Q')-\hiTr(A;Q)=-\rTr\left(A\log(Q'/Q)\right).
\end{equation*}
\end{lemma}

\begin{proof} In view of the cocycle condition it suffices to compute the
difference when $Q$ satisfies the assumption of
Corollary~\ref{rbm.250}. Since $\hiTr(A;Q')$ is the regular value at
$\tau=0$ of the analytic continuation of $\iTr(A(Q')^{-\tau})$ the
difference is the value at $\tau=0$ of $\tau\iTr(AB(\tau)Q^{-\tau})$ where 
\begin{equation*}
B(\tau)=\tau^{-1}\left((Q')^{-\tau}Q^{\tau}-\Id\right)
\end{equation*}
is regular at $\tau=0.$ Thus $B(0)=-\log(Q'/Q)$ and the difference is
$\rTr(AB(0)).$
\end{proof}

\section{Exterior derivations\label{Sect.ED}}

If $D$ is a derivation on an algebra $\mathcal A$ then the map 
\begin{equation*}
i_D(a_0\otimes a_1\otimes a_2 \otimes \cdots\otimes a_p)=a_0 D(a_1)\otimes 
a_2 \otimes \cdots \otimes a_p
\end{equation*}
extends to a map on the Hochschild chain spaces, $i_D:\HoCh p{\mathcal
A}\longrightarrow \HoCh{p-1}{\mathcal A},$ which anticommutes with the
Hochschild differential, $bi_D+i_Db=0.$ It therefore induces a morphism on
Hochschild homology, which we shall also denote as 
\begin{equation*}
i_D:\Hd_p({\mathcal A})\longrightarrow \Hd_{p-1}({\mathcal A}).
\end{equation*}
If $D$ is an inner derivation, given by commutation with an element of
$\mathcal A,$ then this map is trivial.

For the algebras considered here there are two natural exterior
derivations. These are both given formally by commutation with the logarithm of
an operator
\begin{gather}
D_{\log Q}:{\mathcal A}\ni A\longmapsto [\log
Q,A]=\frac{d}{d\tau}\left(Q^{\tau}AQ^{-\tau}-A\right)\big|_{\tau=0}\Mand\\
D_{\log x}:{\mathcal A}\ni A\longmapsto [\log
x,A]=\frac{d}{dz}\left(x^{z}Ax^{-z}-A\right)\big|_{z=0}.
\label{rbm.256}\end{gather}
Here $Q$ is a positive element as considered in Section~\ref{Sect.RF} and
$x$ is a boundary defining function. Since these exterior derivations are
well defined modulo interior derivations, the resulting maps on homology 
do not depend on the choice of $Q$ or $x.$

\section{Explicit cocycles\label{Sect.EC}}

To obtain explicit representations for Hochschild cocycles we begin with
the simple but important case of traces.

\begin{lemma}\label{rbm.258} The functional $\rTr$ spans the 
groups $\Hd^0(\alg1)$ and
$\Hd^0(\alg3),$ $\dTr$ spans $\Hd^0(\idl2)$ and $\iTr$ spans $\Hd^0(\idl4);.$
\end{lemma}

\begin{proof}
The homology groups in dimension $0$ are known to be singly generated
from the computations in Section~\ref{Sect.HH} and the results of
Section~\ref{Sect.RF}. Since the traces $\rTr,\dTr$ and $\iTr$ are all 
non zero, the result follows.
\end{proof}

The arguments in Section~\ref{Sect.HH} can also be used to compute the
Hochschild cohomology groups of the various algebras, as the homology groups
of the corresponding spaces. Using the canonical 
orientations coming from the symplectic structure and the counterclockwise
orientation of the circles, we shall replace homology by cohomology using
Lefschetz and Poincar\'{e} duality. Recalling that we always use complex
coefficients, we therefore obtain isomorphisms
\begin{gather}
\label{rbm.204}
\Phi_1:H^{*}(\cS_{\pa X}X\times\SSs) \oplus
H^{*}_{\rel}(\cS M\times\SSs) 
\simeq \Hd^*(\alg1)\\
\label{rbm.202}
\Phi_2:H^{*}(\cS X\times\SSs) \simeq \Hd^*(\idl2)\\
\label{rbm.203}
\Phi_3:H^{*}(\cS_{\pa X}X\times\SSs\times\SSb)\simeq\Hd^*(\alg3)\Mand\\
\label{rbm.324}
\Phi_4:H^*(\RR \times \SSb) \simeq H^*(\SSb) \simeq \Hd^*(\idl4).
\end{gather}

Denote by $d\theta_r$ the class of $\imath r^{-1}dr$ in $H^*(\SSs)$ and by
$d\theta_x$ the class of $\imath x^{-1}dx$ in $H^*(\SSb).$

\begin{proposition}\label{rbm.257} Under the isomorphism $\Phi_i,$
the action of the morphism $i_{\log Q}$ is transformed to exterior product with
$d\theta_r$ for $i=1,2,3$ and the action of the morphism
$i_{\log x}$ is transformed to exterior product with $d\theta_x$
for $i=3,4:$
\begin{gather*}
\Phi_i(d\theta_r \xi)=i_{\log Q}\Phi_i(\xi), \Mand
\Phi(d\theta_x \xi)=i_{\log x}\Phi(\xi).
\end{gather*}
\end{proposition}

\begin{proof}The isomorphisms $\Phi_i$ are all defined through variants of
the Hochschild-Kostant-Rosenberg map on the symbol algebra. Consider for
example $i=1.$
Denote by $A$ the operator of antisymmetrization in the
variables $f_1,\ldots,f_p$
$$
A(f_0 \otimes f_1 \otimes \ldots \otimes f_p) = (p!)^{-1} \sum_\sigma
\epsilon(\sigma) f_0 \otimes f_{\sigma(1)} \otimes \ldots \otimes f_{\sigma(p)}
$$
where $\epsilon(\sigma)$ is the sign of an element $\sigma$ in the symmetric 
group. Explicitly, if $[\eta]$ is any $p$-cohomology class then the relation 
\begin{equation}
\label{rbm.205}
\Phi([\eta]) \bigl(A(f_0 \otimes f_1 \otimes \ldots \otimes f_p)\bigr)
=\int\limits_{S^*M \times \SSs} *
(\chi_0(f_0 \otimes f_1 \otimes \ldots \otimes f_p))\eta
\end{equation}
for $f_0 \otimes f_1 \otimes \ldots 
\otimes f_p \in F_{n-p}(\alg1)^{\otimes(p+1)}$    
completely determines $\Phi([\eta]).$  We can assume that the principal
symbol of $Q$ is $r.$ The relation between the contraction
$i_{\log Q}$ on the Hochschild complex and the Hochschild-Kostant-Rosenberg
map $\chi_0,$ is expressed  by 
$$
\chi_0\circ i_{\log Q} \circ A (f_0 \otimes f_1 \otimes \ldots \otimes f_p)
= \imath i_{H_{\log r}}\chi_0
(f_0 \otimes f_1 \otimes \ldots \otimes f_p)
$$
where the on the right hand side $i_{H_{\log r}}$ is the contraction with
the Hamiltonian vector field of $\log r.$ The equation
\eqref{rbm.205} proves the result in this case. The other relations follows 
in the same way.
\end{proof}

It follows that $\Hd^1(\idl4) \simeq \CC$ is generated by $\beta = i_{\log
x} \iTr .$ If $B\sim \sum_{-\infty}^m x^{-j}B_j$ and $C \sim
\sum_{-\infty}^m x^{-j}C_j$ then 
\begin{equation}
\beta(B,C)=
\sum\limits_{j,k}\int_{-\infty}^{\infty}\int_{\pa X^2}
(-1)^k k^{-1}B_{-j}(t,y,y')C_{j-k}(t,y',y)t^kdt. 
\label{rbm.123}\end{equation}

We shall now briefly describe the construction of explicit Hochschild cocycles 
corresponding to other forms. For any element
of $(\cPwf{-\ZZ}0 X)^{\otimes(p+1)},$ of the form
$b=f_0 \otimes f_1 \otimes \ldots \otimes f_p$ consider the multilinear map
$\tau_b=\dTr\circ i_b$ 
$$
\tau_b(g_0,\ldots,g_n)=\dTr(f_0g_0[f_1,g_1] [f_2,g_2] \ldots [f_p,g_p]).
$$
This operation extends to any element $b\in\HoCh p{\cPwf{-\ZZ}0 X}.$
Let $\eta$ be a closed form representing a class $[\eta] \in H^{p}(\cS X).$
We can choose the element $b$ so that $\chi_0(b) =\eta,$ where $\chi_0$ is 
the Hochschild-Kostant-Rosenberg map applied to the principal
symbols, and such that $\tau_b$ is a Hochschild cocycle. 
It follows then from the spectral sequence computation in the proof of 
Proposition~\ref{rbm.116} that $\tau_b$ represents the class of $\eta,$
$\tau_b = \Phi_2([\eta]).$ (One can recognize here the action of the
Gerstenhaber algebra $\Hd^*(A,A)$ on the trace $\dTr.$)

The description of the cocycles for $\alg1$ corresponding to classes
$[\eta] \in H^{*}_{\rel}(\cS M\times\SSs)$ is also obtained as
multilinear forms $\tau_b$ or $\tau_b \circ i_{\log Q},$ defined using the
residue trace and tensors $b$ with components vanishing to infinite order
at the boundary. In order to obtain explicit formulae for classes in
$H^{*}(\cS_{\pa X} X\times\SSs) \subset \Hd_*(\alg1),$ it is enough to do so
for $\alg3,$ in view of the previous proposition.

Similar results hold for $\alg3.$ A cocycle $\varphi$ represents the class
of the form $\eta,$ then $\varphi \circ i_{\log x}$ represents the class of
$d\theta_x \eta,$ and $\varphi \circ i_{\log Q}$ represents the class of
$d\theta_r \eta.$ Also the classes not containing any of $d\theta_x$ or
$d\theta_r$ can be obtained as multilinear functionals of the form $\tau_b,$
constructed using invariant operators and the boundary residue $\rTr.$

\section{The homology of ${\mathcal A}_{\pa}$ and $\ialg$\label{Sect.THOA}}

In this section we shall complete the computation of the Hochschild
homology groups announced in the introduction. The main tool 
will be to compute the boundary morphisms in the exact sequence in 
Hochschild homology for an $H$-unital algebra. These results will
be used in the proof of the index formula in the
sections \ref{Sect.BI}, \ref{Sect.IF}. Using also results from 
\cite{Nistor2} we shall give a cohomological interpretation 
of the index formula. 

Next we wish to consider the `Laurent series in finite order indicial
operators' i.e. 
\begin{equation*}
\alg5=x^{-\ZZ}\cPf{\ZZ}X/x^{\ZZ}\cPf{\ZZ}X.
\end{equation*}
Clearly $\idl4\subset\alg5$ is an ideal so consider the short exact sequence
$$
0\longrightarrow\idl4\longrightarrow\alg5\longrightarrow\alg3\longrightarrow0.
$$
Since the structure of this algebra is somewhat more involved we do not
attempt to compute its Hochschild homology directly.
The H-unitality of $\idl4$ means that this gives rise to a long exact
sequence in Hochschild homology which, following Proposition~\ref{rbm.119},
breaks up into 
\begin{gather}
{\begin{split}
0\leftarrow\Hd_0(\alg3)\longleftarrow\Hd_0(\alg5)\longleftarrow\Hd_0(\idl4)
\stackrel{\pa}{\longleftarrow} \Hd_1(\alg3)\longleftarrow\\ 
\Hd_1(\alg5) \longleftarrow\Hd_1(\idl4)
\stackrel{\pa}{\longleftarrow}
\Hd_2(\alg3)\longleftarrow \Hd_2(\alg5)
\longleftarrow 0,\end{split}}\label{rbm.224}\\
\Hd_k(\alg3)\simeq\Hd_k(\alg5),\ k\ge3.\label{rbm.225}
\end{gather}

We shall need the following lemma.

\begin{lemma}
\label{rbm.325}
The following commutation relations hold
\begin{gather*}
\hiTr([A,B])=\rTr(A[\log Q, B])=(i_{\log Q}\rTr)(A,B),\\
\hdTr([A,B])=-\rTr(A[\log x, B])=-(i_{\log r}\rTr)(A,B).
\end{gather*}
\end{lemma}

\begin{proof}By definition, $\hiTr(C)$ is the regular value at
$\tau=0$ of the analytic continuation of $\iTr(CQ^{-\tau}).$ Now for $C=[A,B],$
using the fact that $\iTr$ is a trace, 
\begin{equation*}
\iTr([A,B]Q^{-\tau})=\iTr(A[B,Q^{-\tau}])=\tau\iTr(AM(\tau)Q^{-\tau})
\end{equation*}
where $M(\tau)=\tau^{-1}(B-Q^{-\tau}BQ^{\tau})$ is entire. Since
$M(0)=[\log Q,B]$ it follows from Corollary~\ref{rbm.250} that 
\begin{equation}
\hiTr([A,B])=\rTr(A[\log Q,B])=(i_{\log Q}\rTr)(A,B)
\label{rbm.259}\end{equation}
since $i_{\log Q}(A\otimes B)=A[\log Q,B],$ by definition.

The second relation is proved in a similar manner.
\end{proof}

Consider the dual exact sequence to \eqref{rbm.224}. The boundary map 
\begin{equation}
\label{rbm.323}  
\pa:\Hd^i(\idl4)\longrightarrow \Hd^{i+1}(\alg3)
\end{equation}
is defined explicitly, for $i=0,$ as follows. First the functional $\iTr$ on $\idl4$
should be extended continuously to $\alg5.$ This is already done in
Section~\ref{Sect.RF}, with the result denoted $\hiTr.$ Then the image
$[f]=\pa[\iTr]$ is given by evaluation on the boundary, \ie\ the commutator 
\begin{equation*}
f(A,B)=(b \hiTr)(A,B)=\hiTr([A,B]),\ A,B\in\falg.
\end{equation*}
Here the functional is defined on $\alg5$ since $\hiTr([A,B])=0$ if either
element is in $\idl4$ (otherwise the H-unitality would need to be used
explicitly). By the above lemma $f=i_{\log Q}\rTr.$ 

In particular this shows that $\Phi_3^{-1}(\pa[\iTr])$ is a 
nonzero multiple of $d\theta_r.$
Extending the argument a little gives

\begin{proposition}\label{rbm.260} The boundary maps in \eqref{rbm.224} are
both injective; their duals, \eqref{rbm.323}, are given by
\begin{gather*}
\pa [\iTr] = [i_{\log Q} \rTr]\\
\pa [i_{\log x}\iTr] = -[i_{\log x}i_{\log Q} \rTr].
\end{gather*}

\end{proposition}

\begin{proof} The preceding computation, on the dual complex, shows this
for $i=0.$ Now, apply the morphism $i_{\log x}$ to $\iTr.$ By
Proposition~\ref{rbm.257} we know that the image spans $\Hd^1(\idl4).$

In order to compute $\pa [i_{\log x} \iTr]$
we proceed as above. The functional  $g=i_{\log x} \hiTr$ is an
extension of $g=i_{\log x} \iTr,$ hence 
\begin{gather*}
\pa [i_{\log x} \iTr] = [b i_{\log x} \hiTr] = - [i_{\log x}b \hiTr]
= -[i_{\log x} i_{\log Q}\hiTr].
\end{gather*}
It
follows that the second boundary map is also injective.
\end{proof}

We can now complete the description of the homology of $\alg5.$

\begin{proposition} The Hochschild homology groups of $\alg5$
are
\begin{equation}
\Hd_i(\alg5)\simeq\begin{cases} \ker \pa:\Hd_i(\alg3) \to \Hd_{i-1}(\idl4)
& \Mfor  i=1,2, \\ \Hd_i(\alg3) & \Motwip \end{cases}
\label{rbm.262}\end{equation}
\end{proposition}

\begin{proof}
This is  an immediate consequence of the previous result applied to
the exact sequence \eqref{rbm.224}.
\end{proof}

Combining the results above we can also describe the homology of the
the `most important' algebra
\begin{equation}
\ialg=x^{-\ZZ}\cPf{\ZZ}X/x^{\ZZ}\cPf{-\ZZ}X.
\label{rbm.124}\end{equation}
The exact sequence
\begin{gather}
0\longrightarrow \idl4\longrightarrow\ialg\longrightarrow\alg1\longrightarrow 0
\label{rbm.175}
\end{gather}
determines a long exact sequences in homology.

\begin{proposition}\label{rbm.120} The connecting morphisms 
$\pa:\Hd^i(\idl4) \to \Hd^{i+1}(\alg1)$ 
are given on generators by $\pa[\iTr]=[i_{\log Q} \rTr]$ 
and $\pa[i_{\log x} \iTr] =0;$ it follows that
\begin{gather*}
\Hd^k(\ialg)\simeq \begin{cases}
\CC \oplus H^*_{\rel}(\cS X\times\SSs)\oplus 
H^*(\cS_{\pa X}X\times\SSs)/\CC                                       
&\Mfor k=1\\
H^*_{\rel}(\cS X\times\SSs)\oplus
H^*(\cS_{\pa X}X\times\SSs)
&\Motwip\end{cases}
\end{gather*}
\end{proposition}

\begin{proof}
The naturality of the connecting morphism and the previous proposition
together reduce the proof to the computation of the morphism 
$$p^*:\Hd^*(\alg3) \to \Hd^*(\alg1)$$
induced by the natural map $p:\alg1 \to \alg3.$ 
Using the isomorphisms in the
equations  \eqref{rbm.203} and \eqref{rbm.204}, we deduce from the definition
of the Lefschetz isomorphism that the composition 
\begin{multline*}\Phi_1^{-1} \circ p^* \circ \Phi_3:
H^{k}(\cS_{\pa X} X \times\SSs) \oplus 
H^{k}(\cS_{\pa X} X \times\SSs) d\theta_x =\\
=H^{k}(\cS_{\pa X} X \times\SSs \times \SSb) \to
H^{k}(\cS_{\pa X} X\times\SSs) \oplus H^{k}_{\rel}(\cS M\times\SSs)
\end{multline*}
preserves the direct sum decomposition,  
$\Phi_1^{-1} \circ p^* \circ\Phi_3 =
i^* \oplus (-1)^q\pa,$ where $i$ is the inclusion $i:\cS_{\pa X} X
\times \SSs \to \cS X\times \SSs$ and 
$\pa$ is the connection morphism in the cohomology exact sequence.
The form $d\theta_\sigma$ has an extension to the whole
manifold $\cS X$ and hence $\pa d \theta_\sigma=0.$
Thus we find that
$$
\pa [i_{\log x} \iTr]=p^*\pa [i_{\log x} \iTr]=-p^*[i_{\log x} i_{\log Q} 
\rTr]=-(-1)^q\Phi_1(\pa d\theta_\sigma)=0.
$$
This completes the proof.
\end{proof}

\section{Morita invariance}

In the following sections we will be concerned with index formul\ae\
involving the cocycles studied in the previous sections.  In order to
obtain operators with nontrivial index we need to consider
pseudodifferential operators acting between sections of two bundles $E_+$
and $E_-,$ with isomorphic pull-backs to the cusp-cosphere bundle $\cS
X.$
 
In this section we shall explain how the previous results extend to the
simple case where $E_+=E_-=E$ over $X.$ We shall replace all the algebras
considered above by the corresponding algebras of operators between
sections of $E;$ for instance $\falg$ is to be replaced by $\cPwf{-\ZZ}\ZZ
{X;E}.$ In general by choosing an embedding of $E$ into a trivial bundle,
$E \subset X \times \RR^N,$ we obtain an idempotent $e$ in $M_n(\CI(X)),$
the projection onto the bundle $E,$ and all our algebras $A$ will be
replaced by the bi-ideals (`corners') $eM_N(A)e.$ This will entail changing
the algebras $x^{-\ZZ}{\mathcal G}(\cT X)$ to $x^{-\ZZ}{\mathcal G}(\cT
X,\hom(E))$ and so on. These new algebras are Morita equivalent to the
original ones and hence have the same Hochschild homology.

\begin{proposition} \label{hop1.346}
If $E$ is a vector bundle over 
a compact manifold with boundary then the
Hochschild homology of $\cPwf{-\ZZ}{\ZZ}{X;E}$
and the analogues of $\ialg,$ $\alg1,$ $\alg2,$ $\alg3,$ $\alg4,$ $\alg5$
and $\ridl$ are of the same dimensions as in the case that $E$ is trivial,
since all these algebras are Morita equivalent to their analogues when $E$
is trivial.
\end{proposition}

\begin{proof} An algebra of the form $eM_N(A)e$ is Morita equivalent
to $A$ if $e=e^2 \in M_N(A)$ is an idempotent not contained in any 
proper two-sided ideal. The rest follows from 
the Morita invariance of Hochschild homology see \cite{Connes2} 
and the references therein.
\end{proof}

>From now on $\ialg,$ $\alg1,$ $\alg2,$ $\alg3,$ $\alg4,$ $\alg5$
and $\ridl$ will denote the new algebras, acting between
sections of $E.$ The results on traces extend trivially to the new
algebras as well.

\section{Boundary index\label{Sect.BI}}

As discussed in the boundaryless case in the Introduction,
we will define an `index class' using the regularization discussed in
Section~\ref{Sect.RF}. In fact $\tZ([A,B];z)$ will in general have a pole
at $z=0,$ so we discuss this functional first. Let 
\begin{equation}
\Bif_{x,Q}(A,B)=\lim\limits_{z\to0}z\tZ([A,B];z),\ A,B\in\falg.
\label{rbm.264}\end{equation}

\begin{proposition}\label{rbm.265} The functional $\Bif_{x,Q}$ vanishes if
either component is in $\idl2$ or $\idl4$ and so defines a
functional on $\alg3$ which is a non-trivial Hochschild cocycle mapping under
$\Phi_3^{-1}$ to a multiple of $-d\theta_x+d\theta_r.$
\end{proposition}

\begin{proof} The existence of $\Bif$ depends on the trace property for
$\rTr.$ Consider the identity
\begin{equation}
\tZ([A,B];z)=\Tr([A,B]x^zQ^{-z})+\ha z^2\Tr([A,B]C(z)x^zQ^{-z})
\label{rbm.268}\end{equation}
where $C(z)=z^{-2}(Q^{-z}x^zQ^zx^{-z}-\Id)$ is an entire family with 
\begin{equation}
C(0)=2[\log Q,\log x].
\label{rbm.267}\end{equation}
Thus the second term in \eqref{rbm.268} cannot affect $\Bif$ and
$\Bif(A,B)$ is the value at $z=0$ of
\begin{multline}
z\Tr([A,B]x^zQ^{-z})=z\Tr(B[x^z,Q^{-z}A])\\
=z^2\Tr(BL_1(z)x^zQ^{-z})+z^2\Tr(L_2(z)Bx^zQ^{-z})
\label{rbm.271}\end{multline} 
where $zL_1(z)=[x^z,A]x^{-z}$ and $zL_2(z)=Q^z[Q^{-z},A].$ Thus both $L_1$
and $L_2$ are entire. Moreover $L_1(0)=[\log x,A]$ and $L_2(0)=-[\log Q,A].$
Thus $\Bif$ can be expressed in term of $\rTr$ as
\begin{equation}
\Bif(A,B)=\rTr(B[\log x-\log Q,A])=(-i_{\log x}\rTr+i_{\log Q}\rTr)(A,B).
\label{rbm.266}\end{equation}
Here we have used the notation of Section~\ref{Sect.ED}. The statements of
the proposition now follow from Proposition~\ref{rbm.257}.
\end{proof}

Notice that the formula for this `boundary index functional' can also be
written 
\begin{equation*}
\Bif(A,B)=(\hiTr+\hdTr)([A,B]).
\end{equation*}
It is therefore the image in $\Hd^1(\alg3)$ of the $0$-cocycle $\iTr+\dTr$
on the ideal $\idl2+\idl4$ under the long exact sequence in Hochschild
cohomology arising from the short exact sequence 
\begin{equation*}
0\longrightarrow \idl2+\idl4\longrightarrow\ialg\longrightarrow
\alg3\longrightarrow 0.
\end{equation*}

\begin{proposition}\label{rbm.173} Suppose $A\in\alg3$ is invertible
then the pairing in Hochschild homology 
\begin{equation*}
\cinv(A)=\Bif(A,A^{-1})
\end{equation*}
is homotopy invariant and vanishes if $A$ arises from an invertible element
in $\ialg.$
\end{proposition}

\noindent The first invertibility condition is equivalent to the existence
of a representative in $\cPwf pm{X;E}$ of the form $x^pB$ where $B\in\cPf
m{X;E}$ is elliptic near the boundary; the last invertibility condition
corresponds to the existence of such a $B$ which is globally elliptic and
has invertible indicial family.

\begin{proof} To see the homotopy invariance of $\cinv(A)$ as a function on
the open set of elliptic elements of $\alg3$ consider a smooth family $A_t$
with derivative $\dot A_t.$ The function $\cinv(A_t)$ is smooth while $A_t$
remains elliptic and differentiation gives
\begin{multline*}
\frac d{dt}\Bif_{x,Q}(A_t,A_t^{-1})=\Bif_{x,Q}(\dot A_t,A_t^{-1})-
\Bif_{x,Q}(A_t,A_t^{-1}\dot A_t A^{-1}_t)\\
=\rTr\big(A^{-1}_t\tD{\dot A_t}-A_t^{-1}\dot A_tA^{-1}_t\tD{A_t}).
\end{multline*}
Here 
\begin{equation}
\tD A=[\log x-\log Q,A]
\label{rbm.275}\end{equation}
and \eqref{rbm.266} has been used. Since
$B\mapsto\tD B$ is a derivation and $\rTr$ is a trace this reduces to
$\rTr(\tD{A^{-1}_t\dot A_t}).$ However  
\begin{equation}
\rTr(\tD{B})=0\ \forall\ B\in\falg
\label{rbm.270}\end{equation}
since $F(z)=z^2\Tr(T(z)x^zQ^{-z})\equiv0$ where
$zT(z)=x^zQ^{-z}BQ^zx^{-z}-B$ is entire so $F(0)=\rTr(T(0))=0$ with
$T(0)=\tD{B}.$

Now if $A\in\alg3$ arises from an invertible element $A'\in\ialg,$ \ie\ a
fully elliptic element of $\falg,$ then it has a
parametrix $B'\in\falg$ modulo $\cPwf{\infty}{-\infty}{X;E}.$ In this case
$\cinv(A)$ vanishes since it is the residue at $z=0$ of 
\begin{equation*}
\tilde Z_{x,Q}([A',B'];z)=\Tr([A',B']x^zQ^{-z})
\end{equation*}
which is entire.
\end{proof}

\section{Index $1$-cocycle\label{Sect.I1C}}

For any pair of elements $A,B\in\falg$ we now consider the bilinear
functional corresponding to the regularized value of $\tZ,$ for the
commutator, at $z=0.$ Thus set
\begin{equation}
\If(A,B)=\lim_{z\to0}\big(\tZ_{x,Q}([A,B];z)-z^{-1}\Bif_{x,Q}(A,B)\big),
\ A,B\in\falg.
\label{rbm.168}\end{equation}
The existence of this limit follows from the definition of $\Bif$ in
\eqref{rbm.264}. 

\begin{lemma}\label{rbm.269} If $A\in\falg$ is fully elliptic in the sense
that it is invertible in $\ialg$ then it defines a Fredholm operator on
$\dCI(X;E)$ and its index is  
\begin{equation}
\Ind(A)=\If(A,B),\ B\in\falg,\ [B]=[A]^{-1}\Min\ialg.
\label{rbm.274}\end{equation}
\end{lemma}

\noindent In this case it follows from Proposition~\ref{rbm.173} that
$\tZ([A,B];z)$ is regular near $z=0.$ The fact that $\Ind(A)$ is the value
at $z=0$ is Fedosov's formula. We prove this lemma and discuss the index
formula which follows from it in the next section. This lemma can also be
viewed as expressing the compatibility between the boundary map in
Hochschild homology and the boundary (or index) map in $K$-theory for a
particular cocycle, namely the Fredholm trace. This compatibility is proved
in general in \cite{Nistor2} and will be used in further extensions of the
index theorem, to families for example. See also \cite{Nistor1}.

\begin{proposition} \label{hop1.348}
The index functional defined by \eqref{rbm.168} descends from $\falg$ to 
a cocycle on $\ialg$ and
can be written
\begin{equation}
\If(A,B)=\ha(\hiTr+\hdTr)(B\tD{A}+\tD{A}B).
\label{rbm.198}\end{equation}
\end{proposition}

\begin{proof} By definition $\If(A,B)$ is the regularized value at $z=0$ of
$\tZ([A,B];z).$ Thus is just 
\begin{equation*}
\If(A,B)=\frac{d}{dz} z\tZ([A,B];z)\big|_{z=0}.
\end{equation*}
The symmetrized form of \eqref{rbm.268} gives  
\begin{multline}
2\tZ([A,B];z)\\
=z\Tr\left((BL_1(z)+L_2(z)B)x^zQ^{-z}\right)+
z\Tr\left((BL_3(z)+L_4(z)B)Q^{-z}x^z\right)\\
\Mwith zL_1(z)=x^zAx^{-z}-A,\ zL_2(z)=A-Q^zAQ^{-z},\\
zL_3(z)=Q^{-z}AQ^z-A\Mand zL_4(z)=A-x^{-z}Ax^z.
\label{rbm.276}\end{multline}
The functions $z^2\Tr(Cx^zQ^{-z})$ and $z^2\Tr(CQ^{-z}x^z)$ are both
holomorphic near, and are equal to second order at, $z=0$ with common value
$\rTr(C)$ and derivative $(\hiTr+\hdTr)(C).$ It follows that
\begin{multline*}
2\If(A,B)=(\hiTr+\hdTr)\left(BL_1(0)+L_2(0)B+BL_3(0)+L_4(0)B\right)\\
+\rTr\left(BL'_1(0)+L'_2(0)B+BL'_3(0)+L'_4(0)B\right),
\end{multline*}
Now, $L_1(0)=L_4(0)=[\log x,A],$ $L_2(0)=L_3(0)=-[\log Q,A],$
$L'_1(0)+L'_4(0)=0$ and $L'_2(0)+L'_3(0)=0.$ Since $\rTr$ is a trace
functional, \eqref{rbm.198} follows. This functional certainly vanishes if
either factor is in $\ridl$ so it descends to $\ialg.$
\end{proof}

This computation of the index class can be interpreted in terms of the previous
computations of the Hochschild homology groups. Consider the class of the
functional $\If$ in $\Hd^1(\ialg).$ The restrictions of $\If$ to $\idl2$
and $\idl4$ are given by extension of the Atiyah-Singer cocycle and the
Toeplitz-index cocyle, respectively. These cocycles correspond in our
computations to $i_{\log Q} \dTr$ and, respectively, $-i_{\log x}
\iTr.$ Their sum extends to the whole algebra if and only if the compatibility
condition 
$$\pa i_{\log Q} \dTr = \pa i_{\log x} \iTr
$$
in $\Hd^2(\alg3)$ is satisfied. This checks with the computations of 
lemma~\ref{rbm.325}
and Proposition~\ref{rbm.260} which give that their common value
as the class of the cocycle $i_{\log x} i_{\log Q} \rTr.$ 

Put in an other way, $\If$ gives an extension of $-i_{\log x} \iTr$
from $\idl4$ to $\falg$ and this explains the vanishing of 
$\pa [i_{\log x} \iTr]$ in \ref{rbm.120}.

\section{Index formula\label{Sect.IF}}

\begin{proof}[Proof of Lemma~\ref{rbm.269}] If $A\in\falg$ is fully elliptic
then it has a generalized inverse $B$ which satisfies 
\begin{equation*}
AB-\Id=-\Pi_1,\ BA-\Id=-\Pi_0
\end{equation*}
where $\Pi_0$ and $\Pi_1$ are projections onto the null space and a
complement to the range respectively. Thus $AB-\Id$ and $BA-\Id$ are both
elements of $\ridl$ and 
\begin{multline*}
\Ind(A)=\Tr(\Pi_0)-\Tr(\Pi_1)=-\Tr(BA-\Id)+\Tr(AB-\Id)\\
=-\tZ((BA-\Id);0)+\tZ((AB-\Id);0)=\If(A,B).
\end{multline*}
Now, we know that $\If(A,B)$ only depends on the classes $[A]$ and
$[B]=[A]^{-1}$ in $\ialg,$ so the lemma follows.
\end{proof}

To investigate the index formula \eqref{rbm.274} further we name various of
the component functionals. For any invertible elliptic element $A\in\ialg$ set
\begin{gather}
\etab(A)=-\hiTr(A^{-1}[\log x,A]+[\log x,A]A^{-1})\\
\ASb(A)=\ha\hdTr(A^{-1}[\log Q,A]+[\log Q,A]A^{-1}))\\
\iF(A)=\ha\hiTr(A^{-1}[\log Q,A]+[\log Q,A]A^{-1})\Mand\\
\sF(A)=\ha\hdTr(A^{-1}[\log x,A]+[\log x,A]A^{-1}).
\label{rbm.277}\end{gather}
In fact $\etab(A)$ and $\iF(A)$ are defined for all invertible elements of
$\alg5$ and $\ASb(A)$ and $\sF(A)$ are defined for all invertible elements
of $\alg1.$

With this change of notation we grace the index formula with the name

\begin{Itheorem} \label{rbm.283} For any fully elliptic
element of $\falg$ 
\begin{equation}
\Ind(A)=\ASb(A)-\ha\etab(A)+\iF(A)+\sF(A).
\label{rbm.278}\end{equation}
\end{Itheorem}

\noindent
Here, the notation is supposed to indicate that $\ASb(A)$ is a
generalization of the Atiyah-Singer integrand involving
the $\hat{A}$-genus (or, in this case, the Todd genus) and $\etab(A)$ is a
generalization of the `eta' invariant. Formally it is certainly the case that
$\ASb(A)$ only depends on the (full) symbol of $A$ whereas $\etab(A)$ only
depends on the (full) indicial family of $A,$ \ie\ the respective images in
$\alg1$ and $\alg5.$ To proceed further we make successive simplifying
assumptions on $A.$

\begin{proposition}\label{rbm.279} If $A\in\cPf m{X;E}$ is fully elliptic then
\begin{equation}
\etab(A)=\eta(A_0)\Mwhere A_0=\In(A)
\label{rbm.280}\end{equation}
is the indicial family of $A$ and the `eta' invariant on the right is as
defined in \cite{Melrose46}.
\end{proposition}

The indicial family $A_0=\In(A)$ is defined as the class of $A\in\cPf m{X;E}$ in
the quotient $\cPf m{X;E}/x\cPf m{X;E}.$  If $A\in\falg$ is fully elliptic
then for some $p$ $x^pA\in\cPf m{X;E}$ is elliptic with $A_0$
invertible. Furthermore $A$ and $x^pA$ have the same index. Thus the
assumption of this Proposition involves no essential restriction as far as
the index is concerned.

\begin{proof} If $A\in\cPf m{X;E}$ is fully elliptic then it has a parametrix
$B\in\cPf{-m}{X;E}$ and the indicial families are inverses of each other,
$\In(B)=B_0,$ $B_0=A_0^{-1}.$ For the commutator  
\begin{equation}
[\log x,A]\in x\cPf {m-1}{X;E},\ 
\In(x^{-1}[\log x,A])=-[t,A_0]
\label{rbm.281}\end{equation}
where $A_0\mapsto[t,A_0]$ is the exterior derivation on the indicial
algebra (which is the suspended algebra for $\pa X)$ given by commutation
with the linear function $t.$ Formula \eqref{rbm.195} from
Lemma~\ref{rbm.226} then shows that 
\begin{equation*}
\hiTr(B[\log x,A])=-\bTr(A_0^{-1}[t,A_0])=-\bTr([t,A_0]A_0^{-1})
=\hiTr([\log x,A]B)
\end{equation*}
since $\bTr$ is a trace functional on the full suspended algebra. Thus 
\begin{equation*}
\etab(A)=-\bTr(A_0^{-1}[t,A_0])
\end{equation*}
which is the definition from \cite{Melrose46}.
\end{proof}

\begin{corollary}\label{rbm.282} If $\eth$ is an admissible Dirac operator
on $X$ then 
\begin{equation}
\etab(\eth)=\eta(\eth_0)
\label{rbm.220}\end{equation}
is the Atiyah-Patodi-Singer `eta' invariant of the boundary Dirac operator.
\end{corollary}

\begin{proof} In \cite{Melrose46} the `eta' invariant on the right of
\eqref{rbm.220} is identified with the original definition from
\cite{Atiyah-Patodi-Singer1} for admissible Dirac operators.
\end{proof}

A further simplifying assumption, which can always be arranged by
deformation, is that the operator $A$ is translation invariant near the
boundary.

\begin{lemma}\label{rbm.284} If $A\in\cPf m{X;E}$ is fully elliptic and
translation invariant with respect to some normal fibration, with normal
variable $x'$ near the
boundary then the boundary defining function $x$ can be chosen of the form
$f(x')$ so that $dx$ is supported in the collar neighborhood and if this is
done then
\begin{equation}
\sF(A)=0.
\label{rbm.285}\end{equation}
\end{lemma}

\begin{proof} Certainly taking $x=f(x')$ it can be arranged to be constant
outside any preassigned neighborhood of the boundary. If $A$ is
translation invariant in a neighborhood containing the support of $dx$ then
the Taylor series of the commutator takes the form 
\begin{equation*}
[\log x,A]\simeq\sum\limits_{j=1}^{\infty}\left(\frac{d}{dx'}\right)^jf(x')
R_j
\end{equation*}
with $R_j$ translation invariant and of order at most $m-j.$ The composites
$B[\log x,A]$ and $[\log x,A]B$ therefore have similar expansions with
different coefficients $R_j.$ The explicit formula given in
Proposition~\ref{rbm.233} for the functional $\hdTr,$ as the Hadamard
regularization of Wodzicki's residue trace, shows that all terms with $j>1$
integrate to zero (and only terms with $j\le\dim X$ can be
non-trivial). Furthermore the term with $j=1$ integrates exactly to be
logarithmic, so is removed by the regularization.
\end{proof}

\begin{SItheorem}\label{rbm.287} If $A\in\cPf0{X;E}$ is fully elliptic and
translation invariant with respect to a normal fibration near the boundary
then, for a choice of boundary defining function as in Lemma~\ref{rbm.284},
\begin{equation}
\Ind(A)=\ASb(A)-\ha\eta(A_0)+\iF(A)
\label{rbm.288}\end{equation}
where $\eta(A_0)$ is the `eta' functional of \cite{Melrose46}.
\end{SItheorem}

\noindent As already noted above the assumptions here do not limit the
applicability of the index formula since they can be arranged by deformation.

To simplify the remaining two terms in the index formula appreciably we
need to make a generally non-trivial assumption on the indicial family.

\begin{lemma}\label{rbm.289} If $A \in \cPf m{X;E}$ has invertible indicial family
and is normal to first order at the boundary in the sense that for some
inner product
$$
[A^*,A]\in\cPwf2{2m-1}{X;E}
$$
then $Q$ can be so chosen that
$\iF(A)=0$ and, assuming that $A$ is fully elliptic, the integral for
$\ASb(A)$ requires no regularization.
\end{lemma}

\begin{proof} If $A$ is fully elliptic and of positive order then we can
choose 
\begin{equation*}
Q=(A^*A+1)^{1/2m}. 
\end{equation*}
If it is only elliptic near the boundary, or is
of negative order, then an appropriate operator which is positive in the
interior and vanishes near the boundary can be used in place of $\Id$ so
that $Q$ is positive and of order $1.$ In any case $Q$ can be chosen so
that $[A,\log Q]\in\cPwf2m{X;E}.$ Since $\hiTr$ only depends on the coefficient
of $x$ in the Laurent series of its argument it vanishes on $A^{-1}[\log
Q,A]$ and $[\log Q,A]A^{-1}$ in this case. Similarly the regularization of
the integrand in the definition of $\ASb(A)$ involves only the terms which
do not vanish quadratically.
\end{proof}

Combining these various result we arrive at an index formula which takes
the same form as that of Atiyah, Patodi and Singer obtain for Fredholm
Dirac operators.

\begin{RItheorem}\label{rbm.290} If $A\in\cPf m{X;E}$ is fully elliptic,
translation invariant near the boundary and has normal indicial family
$A_0$ then
\begin{equation}
\Ind(A)=\ASb(A)-\ha\eta(A_0)
\label{rbm.291}\end{equation}
where the Atiyah-Singer integral is given by a convergent integral in terms
of Wodzicki's residue trace  
\begin{equation}
\ASb(A)=\RTr(A^{-1}[\log Q,A])
\label{rbm.293}\end{equation}
for any $Q\in\cPf1{X;E}$ which is positive, translation invariant near the
boundary and commutes there with $A$ and where $\eta(A_0)$ is the eta
invariant on the suspended algebra of the boundary as defined in
\cite{Melrose46}.
\end{RItheorem}

\section{Cyclic homology\label{Sect.CH}}

For the algebras that we have considered the cyclic homology
can be computed directly from the Hochschild homology; we proceed to
explain how to do this.

Recall \cite{Connes2,Loday-Quillen1} that the cyclic homology of a 
unital algebra $A,$ denoted  $\Hc_n(A),$ is the
homology of the cyclic complex $({\mathcal C}(A),b+B)$ where 
\begin{equation}
{\mathcal C}(A)_n=\bigoplus_{k\geq 0}A\otimes(A/\Bbb C 1)^{\otimes n-2k}.
\end{equation}
Completed (\ie\ projective) tensor products are to be used if
topological algebras are considered, as it is the case here. The
differential $B$ is defined by 
\begin{equation}
B(\Tt)=s\sum_{k=0}^{n} t^k(\Tt).
\end{equation}
Here we have used the notation of \cite{Connes2}, that 
$s(\Tt)=1\otimes \Tt$ and 
$t(\Tt)=(-1)^n a_n\otimes a_0\otimes\ldots\otimes a_{n-1}.$

Cyclic and Hochschild homology are related by a long exact sequence due to 
Connes in cohomology \cite{Connes2} and to Loday and Quillen in 
homology \cite{Loday-Quillen1}.
\begin{equation}
\ldots\rightarrow \Hd_{n}(A)\stackrel{I}{\longrightarrow}
\Hc_{n}(A)\stackrel{S}{\longrightarrow}
\Hc_{n-2}(A)\stackrel{B}{\longrightarrow}
\Hd_{n-1}(A)\stackrel{I}{\longrightarrow}\ldots
\end{equation}
$S$ being the periodicity operator.

For the algebras $\alg{1},\idl2,\alg3$ and $\idl4,$ the periodic cyclic
cohomology can be deduced from the Hochschild homology as follows.

\begin{lemma} The operator $B$ in the above long exact sequence vanishes
for the algebras $\alg{1},$ $\idl2,$ $\alg3$ and $\idl4.$ 
\end{lemma}

\begin{proof}
We proceed by induction on $m$ to show that the morphism
$$
B:\Hc_{m-1}(\mathcal A') \to \Hd_m(\mathcal A')
$$
is zero for any $m$ if $\mathcal A'$ is any of the algebras
$\alg{1},$ $\idl2,$ $\alg3$ or $\idl4.$

The statement is trivially true for $m=0.$ We begin with the following remark. 
The proofs of Propositions~\ref{rbm.116}, \ref{rbm.117}, \ref{rbm.118} and
\ref{rbm.219} show that the groups $\Hd_q (\mathcal A')$ are generated by
elements of order $-n+q$ with respect to the order filtration, and that
all cycles of order less that $-n+q$ are boundaries. Here $n$ denotes the 
dimension of the manifold  $M.$ This is a direct consequence of the 
computation of the $E^1$ terms of the spectral sequences associated to the
order filtration. 

Assuming the statement to be true for all values less than a given $m$ we
find, from the Connes' exact sequence, that the groups $\Hc_{m-1}(\mathcal
A')$ are isomorphic to $\Hc_{m-1}(\mathcal A') \simeq \oplus_{k \geq 0}
\Hd_{m-2k-1}(\mathcal A').$ This shows that the groups $\Hc_{m-1}(\mathcal
A')$ are generated by elements of order at most $-n + m -1$ in the order
filtration. It follows that they map to elements of order less than $-n+m$
in $\Hd_m(\mathcal A'),$ and hence they vanish by the remark above.

The last statement is a direct consequence of the vanishing of $B$ 
in the Connes' exact sequence.
\end{proof}

For any exact sequence of algebras we obtain an exact sequence in cyclic
homology \cite{Wodzicki4} completely similar to the exact sequence in
Hochschild homology. Moreover the boundary maps of these two exact
sequences are compatible. We therefore obtain the following result for the
cyclic homology of these algebras.

\begin{proposition} 
\label{hop1.347}
The cyclic homology groups are given by
\begin{equation}
\Hc_m(\mathcal A') =\oplus_{k \geq 0} \Hd_{m-2k}(\mathcal A')
\end{equation}
where $\mathcal A'$ is any of the algebras $\alg{1},\idl2,\alg3,\idl4$ 
or $\alg5.$ Explicitly 
\begin{gather*}
\Hc_j(\alg1)\simeq\oplus_{k \geq 0}\bcoh^{2n-j+2k}(\cS X\times\SSs)
\\
\Hc_j(\idl2)\simeq \oplus_{k \geq 0}H_{\rel}^{2n-j+2k}(\cS X\times\SSs).
\\
\Hc_j(\alg3)\simeq \oplus_{k \geq 0} 
H^{2n-j+2k}(\cS _{\pa X}X\times\SSs\times\SSb)
\\
\Hc_*(\idl4) \simeq \Hd_*(\idl4)\simeq H^{1-j}(\SSb)
\\
\Hc_i(\alg5)\simeq\Hc_i(\alg3)/\CC, \Mif i >0 \;, 
\Hc_0(\alg5)\simeq\Hc_0(\alg3).
\end{gather*}
The periodicity operator $S$ is identified with the canonical projection,
so the periodic cyclic homology groups are obtained by removing the
restriction that $k$ be positive.
\end{proposition}

\begin{proof}
The computations for the first three algebras follow from the computation
of Hochschild homology and the previous lemma. The cyclic homology groups
of $\alg5$ are obtained from the cyclic homology groups of $\alg3$ and
$\idl2$ using the exact sequence in cyclic homology.
\end{proof}

\section{Bundles and grading\label{S.BaG}}

The discussion so far has been concerned with elliptic matrices of
pseudodifferential operators. These results were extended directly to
operators on sections of a particular vector bundle by selecting a
complementary bundle, such that the sum is trivial. However, this still
does not include the case of Dirac operators, since these in general give
operators between non-isomorphic bundles. To overcome this we now briefly
show how the results can be extended to the $\ZZ_2$-graded case.

For a $\ZZ_2$-graded algebra, let $\epsilon$ be the parity operator, so
$\epsilon(a)=(-1)^{\deg a}a$ if $a$ is of pure type. Let $E=E_+\oplus E_-$
be a $\ZZ_2$-graded vector bundle over a compact manifold with boundary
$X.$ Choose a parity-preserving embedding of $E \subset X \times (\CC ^N
\oplus \CC^N)$ into a trivial $\ZZ_2$-graded bundle such that $E=e\bigl(X
\times (\CC ^N \oplus \CC^N)\bigr)$ for an even degree matrix-valued
projecton $e.$ Denote by $\str(a)=\tr(\epsilon a)$ the supertrace
functional on $\hom(E)= e M_{2N}(C^\infty(X)) e.$ The embedding into a
trivial bundle gives to rise a canonical connection $\nabla(\xi)=e d \xi,$ the
Grassman connection on $E.$

Define the 
Hochschild-Kostant-Rosenberg map in this case by 
\begin{equation}
\chi(\Tt)=\str\left(a_0 d a_1\wedge\cdots\wedge d a_n\right)
=\str\left(a_0\nabla a_1\wedge\cdots\wedge\nabla a_n\right)+\ldots
\label{hop1.319}\end{equation}
where the dots represent terms involving the curvature $ede\wedge de$ 
of the Grassman connection. 
 
Reviewing the discussion in
Section~\ref{Sect.FADC} we find

\begin{lemma}\label{hop1.317} If $E=E_+\oplus E_-$ is a $\ZZ_2$-graded
vector bundle over a compact manifold with boundary $X$ the Hochschild
homology, in the graded sense, of the following algebras of sections of the
homomorphism bundle reduce to the corresponding de Rham spaces 
\begin{gather}
\Hd_*(x^{-\ZZ}\CI(X;\hom(E))\simeq x^{-\ZZ}\CI(X;\Lambda^*)\\
\Hd_*(\dCI(X;\hom(E))\simeq \dCI(X;\Lambda^*)
\label{hop1.318}\end{gather}
with identification given by the Hochschild-Kostant-Rosenberg
map $\chi,$ as defined above. This isomorphism does depend on the choice of
a degree preserving embedding  $E \subset X \times (\CC ^N \oplus \CC^N).$
\end{lemma}

\begin{proof} This is proved in essentially the same manner as Proposition~
\ref{hop1.346} using Proposition 1 of \cite{Nistor3}
to first reduce the graded Hochschild homology to the usual Hochschild
homology (of a different algebra). 
\end{proof}

\begin{proposition}\label{hop1.320} If $E=E_+\oplus E_-$ is a $\ZZ_2$
graded vector bundle over a compact manifold with boundary then all the
Hochschild homology groups, 
as $\ZZ_2$-graded algebras, of $\cPwf{-\ZZ}{\ZZ}{X;E}$
and the analogues of $\ialg,$ $\alg1,$ $\alg2,$ $\alg3,$ $\alg4,$ $\alg5$
and $\ridl$ are of the same dimensions as in the case that $E$ is trivial;
moreover the isomorphisms do not depend on the choice of the
embedding of $E$ into a trivial bundle.
\end{proposition}

\begin{proof} These isomorphisms are proved as in the proposition above.
The periodic cyclic homology is independent on the choice of the
embedding since any two embedding are homotopic by a smooth homotopy
(eventually by increasing $N$) and periodic cyclic homology is 
homotopy invariant \cite{Connes2}. Since for  the algebras listed above
Hochschild homology embeds into periodic cyclic homology
(Theorem \ref{hop1.347}) the result follows.
\end{proof}

Suppose that the bundles $E_+$ and $E_-$ are canonically isomorphic in a
neighborhood of the boundary. The embeddings $E_\pm =e_\pm (X \times
\CC^N)$ ($e=e_+ \oplus e_-$) can then be chosen such that their images {\em
coincide} in a neighborhood of the boundary, that is $e_+=e_-$ there.

We now proceed as in section ~\ref{Sect.IF}. We shall  
use the extensions of the traces $\rTr,$ $\iTr,$ $\dTr$ and the functionals
$\hiTr$ and $\hdTr$ to the algebra 
$\cPwf{-\ZZ}{\ZZ}{X;\CC^{2N}}=M_{2N}(\cPwf{-\ZZ}{\ZZ}{X}).$ These extensions
vanish on the off-diagonal terms. Using the embeddings introduced above, we
obtain by restriction functionals on 
$\cPwf{-\ZZ}{\ZZ}{X;E}.$ We shall use the same
notation for these restrictions. Define as before
\begin{equation}
\If(A,B)=\lim_{z\to0}\big(\tZ_{x,Q}([A,B];z)-z^{-1}
\Bif_{x,Q}(A,B)\big),
\ A,B\in \cPwf{-\ZZ}{\ZZ}{X;E}.
\end{equation}
This cocycle descends to a cocycle on $\ialg,$ and
Proposition \ref{hop1.348} extends without change to this case
giving 
\begin{equation}
\label{hop1.349}
\If(A,B)=\ha(\hiTr+\hdTr)(B\tD{A}+\tD{A}B).
\end{equation}
We will also need
a correction term $\nu$ which potentially 
depends on the embeddings, which we now 
proceed to define. Denote by $Q_N=Q \otimes I$ the operator $Q$ extended
to act diagonally on the trivial graded-bundle 
$X \times (\CC^N \oplus \CC^N).$ The functions $Tr(e_\pm Q_N^{-s})$
have no pole at $0$ because $e_+$ and $e_-$ are differential operators
\cite{Wodzicki4}. Then define
\begin{equation}
\nu = Tr(e_+Q_N^{-s})-Tr(e_-Q_N^{-s})\big|_{s=0}
=\str (e Q_{2N}^{-s})\big|_{s=0}
\end{equation}

\begin{theorem}\label{hop1.321} If $E$ is a $\ZZ_2$-graded vector bundle
over a compact manifold with boundary and $A\in\cPwf{-\ZZ}{\ZZ}{X;E}$ is
fully elliptic and odd with respect to the grading then
$[A]\otimes[A]^{-1}$ is a $\ZZ_2$-graded Hochschild cycle and 
\begin{equation}
\Ind(A)=\dim\null(A_+)-\dim\null(A_-)=\If([A],[A]^{-1})+ \nu.
\label{hop1.322}\end{equation}
where $\nu$ is as defined above and the cocycle $\If$ satisfies
\eqref{hop1.349}.
If $A$ is self-adjoint and translation-invariant near the boundary then the
analogue of \eqref{rbm.291} holds.
\end{theorem}

\begin{proof} We consider as before the function $\Tr((e_+ - BA) Q_N^{-s})-
\Tr((e_- - AB) Q_N^{-s})$ which is defined and analytic for large
$s,$ and extends to a holomorphic function at $0$ whose value 
at $0$ is $\Ind(A)$ the index of $A.$ We obtain  
\begin{equation}
\Ind(A)=\Tr((e_+ - BA - e_- + AB) Q_N^{-s})\big|_{s=0}
=\If([A],[A]^{-1})+ \nu
\end{equation}
All the other statements are proved in exactly the same
way as their analogues for $E_+=E_-.$
\end{proof}

Although this theorem involves an assumption on the `boundary triviality'
of the operator $A$ this is satisfied by admissible Dirac operators and in
this sense it is a direct pseudodifferential extension of the
Atiyah-Patodi-Singer index theorem.  This is reflected in our assumption
that $E_+ = E_-$ in a neighborhood of the boundary.

All our previous results on the cohomology classes of the various cocycles
appearing in the index formula extend to the above case ($E_+ \not = E_-$)
with identical proof. The correction term $\nu$ is an interior term, in the
sense that it can be computed as the integral of a local term supported in
the interior of the manifold $X.$ This shows that the identification of the
eta-invariant with one of the terms appearing in our index cocycle
extends to this case as well.

\section{Other algebras\label{Sect.OA}}

The methods employed above can also be used to compute the Hochschild
homology of the algebras $\cPf\ZZ X,$ $\cPf0 X,$ algebras of
complete symbols and the $b$-calculus algebras \cite{Mazzeo-Melrose4}.
However these groups can be infinite dimensional, and hence the relation
with the homology is not readily seen. As an example consider the cusp
operators of order $0,$
$$
{\mathcal A}_6
=\cPwf{-\ZZ}0X/\cPwf{-\ZZ}{-\infty}X.
$$

\begin{proposition} \label{rbm.206} The $E^1_{p,q}$-term of the
spectral sequence associated to the order filtration ${\mathcal A}_6$ is
the same as that of $\alg3$ for $p \le 0$ and vanishes for $p> 0.$ 
The $E^2_{p,q}$-term is unchanged for $p < 0$ and vanishes for $p> 0.$
The $E^2_{0,q}$-terms  are given by 
$$
E^2_{0,q} \simeq
\CI(\cS X;\Lambda^{2n-1-q})\oplus\CI(\cS_{\pa X}X;\Lambda^{2n-q-2})
$$ if $q \not = n$ and
$$
E^2_{0,n} \simeq Z^{n-1}(\cS X) \oplus Z^{n}(\cS X) \oplus
Z^{n-2}(\cS_{\pa X} X) \oplus Z^{n-1}(\cS_{\pa X} X) .$$
The spectral sequence degenerates at $E^2.$
\end{proposition}

\noindent Here have denoted by $Z^k(M)=\ker d,$ the space of closed forms
of degree $k$ on the manifold $M.$

\begin{proof}
The statement about the $E^1$-term follows from the definition.
The differential $d_1: E^1_{p,q} \to E^1_{p-1,q}$ turns out to be the same
if $p \le 0$ and to vanish otherwise. Using the computations
in the proof of proposition \ref{rbm.116}, we see that we need to compute
the homology of the truncated de Rham complexes 
\begin{multline}
0 {\longleftarrow}{\mathcal G}^{q'}(\cS X \times \SSs;\Fo{2n})
{\longleftarrow}{\mathcal G}^{q'}(\cS X \times \SSs;\Fo{2n-1})
{\longleftarrow} \cdots\\
{\longleftarrow}
{\mathcal G}^{q'}(\cS X \times \SSs;\Fo{2n-q})\longleftarrow 0
\end{multline}
where $q'=n-q$ and ${\mathcal G}^k(\cS X\times\SSs;\Fo j)$ denotes the
space of $k$-homogeneous forms. The result then follows.
\end{proof}

\appendix
\section{Cusp calculus\label{App.A}}

We give a definition of the `cusp' calculus on a manifold with boundary, by
use of blow up techniques. This is also done in somewhat greater generality
in \cite{Mazzeo-Melrose4}. Here we describe the construction and refer to
\cite{Mazzeo-Melrose4} for more details. A more traditional
characterization of cusp pseudodifferential operators in terms of standard
pseudodifferential operators on Euclidean space is also given and the basic
properties of the calculus are described.

An algebra of cusp pseudodifferential operators on a compact manifold
with boundary, $X,$ is determined by a choice of trivialization of the
normal bundle to the boundary, which need not be connected. This is
equivalent to the choice of a defining function $\rho$
up to a positive constant multiple near each boundary hypersurface and the
addition of a term vanishing quadratically at the boundary.
We call such a choice a {\em cusp structure} on $X.$

\begin{lemma}\label{rbm.126} Given any two cusp structures on a manifold
with boundary $X$ there is a diffeomorphism of $X,$ connected to the
identity through diffeomorphisms, reducing one to the other.
\end{lemma}

\begin{proof} Let $\rho$ be a defining function inducing one of two given
cusp structures on $X.$ Then, by scaling $\rho$ by a constant as necessary
near each boundary hypersurface it can be assumed that there is a product
decomposition of a neighborhood of the boundary of $X$ of the form 
\begin{equation*}
[0,1]_{\rho}\times Y,\ Y=\pa X.
\end{equation*}
Let $V$ be the normal vector field for this product decomposition, so
satisfying $V\rho =1$ in $\{\rho \le1\}.$ A defining function for the other
cusp structure is necessarily of the form $a\rho,$ with $\CI(X)\ni a>0.$ Only
the restriction of $a$ to the boundary affects the cusp structure, so it
can be assumed that the second defining function is of the form $\rho
'=a\rho$ with $a=1$ near $\rho =1$ and $V\rho'>0$ in
$\{\rho\le1\}=\{\rho'\le1\}.$ It follows that $\rho _t=(1-t)\rho +t\rho '$
is a $1$-parameter family of defining functions with $V(\rho_t)=a_t>0$ in
$\{\rho \le1\}.$  A suitable diffeomorphism can be defined by integration
of the vector field $V_t={a'_t}^{-1}V$ from the fixed surface
$\{\rho_t=1\}$ towards the boundary. 
\end{proof}

The choice of a cusp structure on $X$ defines a Lie algebra of smooth
vector fields by  
\begin{equation}
\Vc(X)=\left\{V\in\CI(X;TX);V\rho\in\rho^2\CI(X)\right\}
\label{rbm.127}\end{equation}
The space of pseudodifferential operators on $X$ described here corresponds
to the {\em microlocalization} of this Lie algebra of vector fields. The
local structure of these vector fields follows by elementary computation.

\begin{lemma}\label{rbm.128} In local coordinates,
$x,y_1,\dots,y_{n-1}$ near a boundary point, with $x$ being equal to an
admissible defining function, an element of $\Vc(X)$ takes the form 
\begin{equation*}
V=x^2\frac{\pa}{\pa x}+\sum\limits_{j=1}^{n-1}b_j\frac{\pa}{\pa y_j}
\end{equation*}
where the coefficients $a$ and $b_j$ are smooth in the coordinate
patch. Conversely a vector field on $X$ is in $\Vc(X)$ if it is of this
form with respect to a covering of the boundary by such coordinate systems.
\end{lemma}

\begin{proof} The last part of the lemma follows from the first part and
the elementary fact, directly from the definition in \eqref{rbm.127}, that
$\Vc(X)$ is a $\CI(X)$ module.
\end{proof}

The local form of the vector fields allows \eqref{rbm.127} to be reversed
in the sense that a defining function $\rho'$ defines the same cusp
structure as $\rho$ if 
\begin{equation*}
\Vc(X)\rho'\in\rho^2\CI(X)
\end{equation*}
where $\Vc(X)$ is defined from the $\rho.$ Said a different way this is

\begin{corollary}\label{rbm.129} The Lie algebra $\Vc(X)$ determines the cusp
structure from which it is defined.
\end{corollary}

The first definition we give is of the `ignorable' operators in the cusp
calculus, these are actually the same for all cusp structures on the given
manifold.

\begin{definition}\label{rbm.130} The space $\igc X$ consists of the
integral operators with kernels in $\dCI(X^2;\pi_R^*\Omega),$ where
$\pi_R:X^2\longrightarrow X$ is projection onto the second factor and
$\Omega$ is the density bundle on $X,$ so 
\begin{equation*}
Af(z)=\int_{X}A(z,z')f(z'),\ f\in\dCI(X)
\end{equation*}
where the same letter is used for operator and kernel.
\end{definition}

\noindent
By Fubini's theorem these operators form an algebra.

The first, and intrinsically global, definition we give of the cusp
calculus uses the properties of real blow-up of p-submanifolds of a
manifold with corners. First consider the space on which kernels of
b-pseudodifferential operators become simple, 
\begin{equation}
\Sp bX2=[X^2;\{H\times H\}_{H\in\bH X}]
\label{rbm.131}\end{equation}
Here $\bH X$ is the set of boundary hypersurfaces of $X.$ The sets of
their products $H\times H$
forms a disjoint collection of p-submanifolds of $X^2,$ so the blow up is
well defined independent of order. As a manifold with corners $\Sp bX2$ has
an additional boundary hypersurface (compared to $X^2)$ for each $H.$ These
faces are each naturally isomorphic (by definition) to the inward-pointing
spherical normal bundle to the corresponding $H\times H\in\bF2{X^2}.$ The
choice of cusp structure on $X$ induces an analogous choice of cusp
structure on $X^2,$ since $\bH{X^2}=\bH X\sqcup\bH X,$ and this provides a
trivialization of the normal bundle to $H\times H.$ Thus, using the given
cusp structure, the extra faces of $\Sp bX2$ are of the form
$[-1,1]\times H^2.$ The `fibre diagonal submanifolds' 
\begin{equation*}
\fD2H=\{0\}\times H^2\subset[-1,1]\times H^2,\ H\in\bH X
\end{equation*}
are therefore fixed by the choice of cusp structure. These p-submanifolds
are disjoint in $\Sp bX2,$ so we may define, without dependence
on the order 
\begin{equation}
\Sp cX2=[\Sp bX2;\{\fD2H\}_{H\in\bH X}],\ \Sp {cb}\beta2:\Sp
cX2\longrightarrow \Sp bX2,
\label{rbm.132}\end{equation}
a manifold with corners defined from $X$ and the choice of cusp structure.

The following basic result is proved in \cite{Mazzeo-Melrose4}.

\begin{lemma}\label{rbm.133} The `stretched projections' $\Sp{c,S}\pi2:\Sp
cX2\longrightarrow X$ for $S=L,R,$ defined as the composites of the
blow-down maps and projection onto the left or right factor 
\begin{equation}
\xymatrix{\Sp cX2\ar@/^2pc/[rrr]^{\Sp{c,S}\pi2}
\ar@/_2pc/[rr]_{\Sp c\beta2}\ar[r]^{\Sp{cb}\beta2}&
\Sp bX2\ar[r]^{\Sp b\beta2}&X^2\ar[r]^{\pi_S}&X},\ S=L\text{ or }R,
\label{rbm.134}\end{equation}
are both b-fibrations.
\end{lemma}

Let $\dCI_{\cff}(\Sp cX2)\subset\CI(\Sp cX2)$ be the subspace of functions
which vanish to infinite order at all boundary hypersurfaces {\em other
than} those, $\cff(H)$ for $H\in\bH X,$ produced in the second stage of blow 
up, i.e. in \eqref{rbm.132}.

\begin{lemma}\label{rbm.135} The closure of the inverse image under $\Sp
c\beta2$ of the diagonal over the interior of $X$ is a p-submanifold,
$\cDiag\subset\Sp cX2,$ meeting the boundary only in the set
$\cff\subset\pa\Sp cX2.$
\end{lemma}

Both the notion of a conormal distribution and the more refined notion of a
polyhomogeneous, or classical, conormal distribution with respect to a
p-submanifold is defined, for example by extension across the boundary. The
space of 1-step polyhomogeneous conormal
distributional sections of order $m$ of a vector bundle $E$ where the
conormality is with respect to a p-submanifold $M\subset X$ of a manifold
with corners $X,$ will be denoted $\pcs mXME.$ Those vanishing to all orders
at boundary hypersurfaces other than those in the subset $S\subset\bH X$
will be denoted $\vpcs mXMES\subset\pcs mXME.$ In fact in terms of the
$\CI(X)$ module structure of $\pcs mXME,$  
\begin{equation*}
\vpcs mXMES=\dCI_S(X)\cdot \pcs mXME.
\end{equation*}

\begin{definition}\label{rbm.136} The space of cusp pseudodifferential
operators, acting on functions, corresponding to a given cusp structure on a
manifold with boundary $X$ is the space of kernels 
\begin{equation}
\begin{aligned}
\cPf mX\equiv&\vpcs m{\Sp cX2}{\cDiag}{\Kd}{\cff}\text{ with}\\
&\Kd=\chd\otimes
(\Sp {c,L}X2)^*(\rho^{\ha}\mhd)\otimes(\Sp {c,R}X2)^*(\chd),
\end{aligned}
\label{rbm.137}\end{equation}
and where $\Sp cX2$ is the blown-up space determined by the cusp
structure.
\end{definition}

\begin{lemma}\label{rbm.139} Let $\mathcal G\subset\CI((\ins X)^2)$ be
such that $\dCI_{\cff}(\Sp cX2)$ is a $\mathcal G$-module with respect to
point wise multiplication, then the space of kernels in \eqref{rbm.137} is a
$\mathcal G$-module over $\CI(X^2);$ it is also invariant under exchange of
the two factors.
\end{lemma}

Since the lift of $\CI(X^2)$ to $\Sp cX2$ has the module structure required
of $\mathcal G$ in this lemma, $\cPf mX$ is a $\CI(X^2)$-module. This
allows the space of kernels acting `from sections of one vector bundle over
$X$ to sections of an other vector bundle over $X$' to be defined as a
tensor product.

\begin{definition}\label{rbm.140} The space of cusp pseudodifferential
operators, acting from sections of any given vector bundle $E$ over $X$ to
sections of any vector bundle $F$ over $X$ and corresponding to a given
cusp structure on the manifold with corners $X$ is the space of kernels 
\begin{equation*}
\cP mX{E,F}=\cPf mX\otimes_{\CI(X^2)}\CI(X^2;\Hom(E,F)) 
\end{equation*}
where $\Hom(E,F)$ is the bundle over $X^2$ with fibre
$\hom(E_{z'},F_{z})\equiv F_{z}\otimes E_{z'}'$  at $(z,z')\in X^2.$
\end{definition}

We shall abbreviate the name to c-pseudodifferential
operator and $\cP mX{E,E}$ to $\cP mXE.$ For any smooth map of compact
manifolds with corners $F:X\longrightarrow Y,$ the push-forward operation
is well defined on supported distributions. If $G$ is any bundle over $Y$
then 
\begin{equation}
F_*:\dCmI(X;F^*G\otimes\Omega)\longrightarrow\dCmI(Y;G\otimes\Omega).
\label{rbm.141}\end{equation}
Here, $\dCmI(X;L),$ for any vector bundle $L,$ is the dual of
$\CI(X;\Omega\otimes L').$ Since the elements of $\cP *X{E,F}$ are smooth
up to all boundaries (in the normal variables to those boundaries), 
\begin{equation}
\cPf *X\subset\dCmI(\Sp cX2;\Kd).
\label{rbm.142}\end{equation}

\begin{proposition}\label{rbm.138} The kernel density bundle on $\Sp cX2$
has a natural identification over $\ins{(\Sp cX2)}$ with $(\Sp
c{\pi_L}2)^*\Omega^{-1}\otimes\Omega$ so \eqref{rbm.141} and \eqref{rbm.142}
give a push-forward map 
\begin{equation*}
\cPf *X\longrightarrow \dCmI(X)
\end{equation*}
which in fact takes values in $\CI(X).$
\end{proposition}

\begin{lemma}\label{rbm.143} Each c-pseudodifferential operator
$A\in\cP*X{E,F}$ defines consistent continuous linear maps
\begin{gather}
A:\CI(X;E)\longrightarrow\CI(X;F),\label{rbm.146}\\
A:\dCI(X;E)\longrightarrow\dCI(X;F),\\
A:\dCmI(X;E)\longrightarrow\dCmI(X;F),\\
A:\CmI(X;E)\longrightarrow\CmI(X;F).
\label{rbm.144}\end{gather}
\end{lemma}

The consistency here is with respect to the inclusions and projections 
\begin{equation}
\xymatrix{&&\dCmI(X;E)\ar@{>>}[dd]\\
\dCI(X;E)\ar@{^{(}->}[r]&\CI(X;E)\ar@{^{(}->}[ur]\ar@{^{(}->}[dr]\\
&&\CmI(X;E).}
\label{rbm.145}\end{equation}

If $\rho_H$ is a defining function for the boundary hypersurface $H\in\bH
X$ then $g=\pi_L^*(\rho_H^{-1})\pi_R^*\rho_H$ is smooth on $X^2$ except 
at the boundary face $H\times X.$ Moreover it is equal to $1$ on the
diagonal and lifts to $\Sp bX2$ to be smooth up to the front face,
i.e. there is again only one hypersurface up to which it is not smooth. The
same holds for the lift to $\Sp cX2,$ and since the singularity is of
finite order $g$ is a multiplier on $\cPf mX.$ Thus
\begin{equation*}
\cPf mX\ni A\longmapsto\rho_H^{-1}\cdot A\cdot\rho_H\in\cPf mX
\end{equation*}
is a morphism. Combined with \eqref{rbm.146} this shows that there is a
unique operator $A_H$ completing the commutative diagrame with vertical
maps given by restriction 
\begin{equation}
\xymatrix{\CI(X)\ar@{>>}[d]\ar^{A}[r]&\CI(X)\ar@{>>}[d]\\
\CI(H)\ar^{A_H}[r]&\CI(H).}
\label{rbm.155}\end{equation}

This can be generalized by considering the more singular function 
\begin{equation}
\mu_H=\frac1{\pi_L^*\rho_H}-\frac1{\pi_R^*\rho_H}.
\label{rbm.153}\end{equation}
On $X^2$ this is singular at both $H\times X$ and $X\times H.$ Written in
terms of $s_H=\pi_L^*\rho_H/\pi_R^*\rho_H$ and
$r_H=\pi_R^*\rho_H+\pi_L^*\rho_H$ it is
\begin{equation}
\mu_H=\frac{1-s_H^2}{r_Hs_H}
\label{rbm.154}\end{equation}
The surface blown up (in this boundary face) to produce $\Sp cX2$ is
$s_H=1,r_H=0.$ It follows that $\mu_H$ is smooth up to the new boundary
face so produced and hence up to all boundary hypersurfaces of $\Sp cX2$
except the three arising from the lifts of $H\times X,$ $X\times X$ and the
front face of $\Sp bX2$ corresponding to $H.$ Since $\mu_H$ is real and
only has a finite order singularity at these faces, it follows that for any
$\zeta\in\RR$ the function $\exp(i\zeta\mu_H)$ is a multiplier on $\cPf
mX.$ Lemma~\ref{rbm.139} then shows that conjugation gives a smooth map
\begin{equation}
\RR\ni\zeta\longmapsto \exp(i\zeta/\rho_H)A\exp(-i\zeta/\rho_H)\in\cPf
mX,\Mfor A\in\cPf mX.
\label{rbm.148}\end{equation}

\begin{definition}\label{rbm.147} The indicial family of an element
$A\in\cPf mX$ at a boundary hypersurface $H\in\bH X,$ fixed by the choice
of a boundary defining function $\rho_H,$ is the family of operators on
$\CI(H)$
\begin{equation*}
\In_H(A,\zeta)=\big(\exp(\frac{i\zeta}{\rho_H})A
\exp(-\frac{i\zeta}{\rho_H})\big)_H
\end{equation*}
where the operator restriction is defined by \eqref{rbm.155}.
\end{definition}

Consider the special case of this construction for the manifold $[-1,1].$ The
interior of this interval can be identified with $\RR$ by radial
compactification. That is, consider the diffeomorphism
\begin{equation}
\SP:\RR\ni z\longmapsto\frac z{(1+|z|^2)^{\ha}}\in(-1,1)\subset[-1,1].
\label{rbm.150}\end{equation}
This `fixes a cusp structure on $[-1,1]$', in fact in dimension one there is a
unique cusp structure.

\begin{proposition}\label{rbm.149} The c-pseudodifferential operators on
$[-1,1]$ are identified by \eqref{rbm.150} with the pseudodifferential
operators on $\RR$ with `polyhomogeneous coefficients', that is  
\begin{equation}
Au(z)=\frac1{2\pi}\int_{\RR}e^{i(z-z')\zeta}a(\ha(z+z'),\zeta)u(z')dz'
\label{rbm.151}\end{equation}
where the amplitude $a$ is determined by the fact that there is an element
$a'\in\CI([-1,1]\times\RR)$ which is a polyhomogeneous symbol in the
usual sense and such that 
\begin{equation}
a=(\SP\times\Id)^*a'.
\label{rbm.152}\end{equation}
\end{proposition}

\noindent The `usual sense' in which $a'(x,\zeta)$ is a polyhomogeneous
symbol is that $a'\in\CI([-1,1]\times\RR)$ and there are functions
$a'_j\in\CI([-1,1]\times\{-1,+1\})$ such at  
\begin{equation*}
a'\sim\sum\limits_{j=m}^{-\infty}|\zeta|^ja_j(x,\frac\zeta{|\zeta|})\Mas|\zeta
|\to\infty,
\end{equation*}
meaning in turn that for each $N$ and all $\alpha ,\beta \in\NN_0,$
\begin{equation*}
\sup_{|\zeta|\ge1}|\zeta|^{N+1+|\beta|}
|D^\alpha _xD^\beta_\zeta\big(a'-\sum\limits_{j=m}^{-N}|\zeta|^ja_j\big)|
<\infty.
\end{equation*}

\begin{proof} Since $\mu_H$ and $r_H,$ give
coordinates near $\cDiag$ the problem is reduced to operators of order
$-\infty$ where it is a matter of simple estimates arising from the
infinite order vanishing.
\end{proof}

This identification in the case of the model manifold $[-1,1]$ suggests a
characterization of c-pseudodifferential operators by localization. For an
open set $U\subset\RR^{n-1}$ consider pseudodifferential operators on
$\RR\times U$ with kernels of the form 
\begin{equation}
A(z,y,z',y')=(2\pi)^{-n}\int_{\RR^n}e^{i(z-z')\zeta+(y-y')\cdot\eta}
a(z,z',y,y',\xi,\eta)d\zeta d\eta 
\label{hop1.335}\end{equation}
where $a=(\SP\times\SP\times\Id)^*a'$ with $a'\in\CI([-1,1]^2\times
U^2\times\RR^n)$ a symbol in the same sense, 
\begin{equation}
a'\sim\sum\limits_{j=m}^{-\infty}|(\zeta,\eta)|^j
a_j(x,x',y,y',\frac{(\zeta,\eta)}{|(\zeta,\eta)|})\Mas|(\zeta,\eta)|\to\infty.
\label{hop1.336}\end{equation}

\begin{proposition}\label{hop1.334} If $(x,y)$ are coordinates in a
coordinate patch $O=[0,1)\times U.$ $U\subset\RR^{n-1},$ near a boundary
point of a compact manifold with boundary, with $x$ an admissible boundary
defining functions for a given cusp structure, then a kernel as in
\eqref{hop1.335} and \eqref{hop1.336} such that $A'=(T\times T)^*A,$
$T(x,y)=(\SP^{-1}(x+1),y)$ which has relatively compact support in $O^2$ in
these local coordinates defines an element of $\cPf mX.$ Conversely any
element of $\cPf mX,$ for $X$ a compact manifold with boundary is a finite
sum of such local operators, plus an element of $\cPf{-\infty}X$ and a
classical pseudodifferential operator of order $m$ with kernel compactly
supported in the interior of $X^2.$
\end{proposition}

Although this proposition reduces cusp pseudodifferential operators to
rather familiar objects on $\RR\times U$ the global representation
discussed above is very convenient. In particular the conormal
distributions to a submanifold (in this case transversal to the boundary)
always have a global representation, modulo smooth functions, in terms of symbols.

\begin{proposition}\label{hop1.337} For any compact manifold with boundary
there is a global quantization map  
\begin{equation}
x^{-\ZZ}S^{\ZZ}(\cT X)\longrightarrow \cPwf{-\ZZ}{\ZZ}X
\label{hop1.340}\end{equation}
which is order-filtered and induces the isomorphism \eqref{rbm.306}.
\end{proposition}

\begin{proof} Choose a Riemann metric on $\Sp cX2$ with respect to which
the boundary hypersurfaces are totally geodesic (this is always possible on
a compact manifold with corners, see for example
\cite{Hassell-Mazzeo-Melrose2}). Let $\phi\in\CI(\Sp cX2)$ be identically
one in a neighborhood of the lifted diagonal and have support in a collar
neighborhood, $O,$ defined by normal geodesic flow with respect to this
metric. Writing the Riemannian identification as 
\begin{equation*}
G^{-1}:\ct X\equiv N\Diag\supset O'\longrightarrow O\subset \Sp cX2
\label{hop1.339}\end{equation*}
the quantization map is defined by Fourier transformation on the fibres of
$\cT X$
\begin{equation}
A_a=\phi G^*(A'_a)\nu,\ A'_a(p,v)=
(2\pi)^{-n}\int_{\cT X}e^{-iv\cdot\xi}a(p,\xi)d\xi.
\label{hop1.338}\end{equation}
Here $\nu$ is a non-vanishing section of the kernel density bundle and
$d\xi$ is the Riemannian volume form on the fibres of $\cT X,$ which is
identified with the conormal bundle to the lifted diagonal by the lifting
of cusp vector fields from the left factor.

That this map induces an isomorphism \eqref{rbm.306} follows from the fact
that the Fourier transform in the normal bundle gives an identification of
conormal distributions at the submanifold, modulo \ci\ functions, with
symbols modulo Schwartz functions.
\end{proof}

Now $\cPwf{-\ZZ}{-\infty}X$ is an ideal, so \eqref{hop1.340} induces a
product on $x^{-\ZZ}\CI(\cS X)[q]]$ which we denote by $\star.$ It is
indeed a star product in the sense of deformation quantization
. That is, 
\begin{equation}
a\star b=\sum\limits_{j=0}^{\infty}P_j(a,b)
\label{hop1.341}\end{equation}
where $P_j$ is a bilinear differential operator of total homogeneity $-j$
satisfying
\begin{equation}
P_0(a,b)=ab,\ P_1(a,b)-P_1(b,a)=\frac1i\{a,b\}.
\label{hop1.342}\end{equation}
Here, $\{,\}$ is the degenerate Poisson bracket on functions on $\cT X$
induced by the singular symplectic structure which arises from its
identification with $T^*X$ over the interior. The form \eqref{hop1.341}
follows directly from the symbol calculus and then the explicit
identifications in \eqref{hop1.342} follow by continuity from the
interior. Notice that in terms of the admissible coordinates $x,y$ near a
boundary point and the induced canonical coordinates $(x,y,\tau,\eta)$ in $\cT X$
\begin{equation}\begin{gathered}
\cT X\ni\alpha =\tau\frac{dx}{x^2}+\eta\cdot dy,\\
\{a,b\}=\pa_{\tau}a\cdot x^2\pa_{x}b-x^2\pa_xa\cdot\pa_{\tau}b+
\pa_\eta a\cdot\pa_yb-\pa_ya\cdot\pa_{\tau}b.
\end{gathered}\label{hop1.343}
\end{equation}

The algebra $\cPwf{-\ZZ}{\ZZ}X$ is also filtered by the boundary order; the
corresponding residual ideal being $\cPwf{\infty}{\ZZ}X.$ The global
representation of the operators allows the quotient to be readily
identified. To do so, consider a compact manifold without boundary, $Y.$
The `suspension' of $Y$ is the product $\RR\times X.$ Using the
stereographic compactification of $\RR$ to $[-1,1]$ via \eqref{rbm.150}
$\RR\times Y$ can is identified with the interior of the compact manifold
with boundary $Y_{\sus}=[-1,1]\times Y.$ The action of $\RR$ on $\RR$ by
translation lifts to a smooth action of $\RR$ on $[-1,1]$ by projective
linear transformations. Consider the subalgebra of $\cPf{\ZZ}{Y_{\sus}}$
consisting of the operators which are invariant under this $\RR$-action, we
shall denote this algebra by $\Psu{\ZZ}Y$ since it is precisely the
`suspended algebra of pseudodifferential operators' considered in
\cite{Melrose46}. The choice of a boundary defining function (admissible for
a given cusp structure) trivializes the normal bundle to the boundary
and so gives an identification 
\begin{equation*}
(\pa X)_{\sus}\equiv N\pa X.
\end{equation*}

\begin{proposition}\label{hop1.344} The choice of a boundary defining
function admissible for a given cusp structure on a compact manifold with
boundary and a normal fibration of the manifold near the boundary induces
an isomorphism 
\begin{equation}
\cPwf{-\ZZ}{\ZZ}X/\cPwf{\infty}{\ZZ}X\equiv \Psu{\ZZ}{\pa X}[[x]
\label{hop1.345}\end{equation}
where the linear variable on the suspension of $\pa X$ can is identified
with $1/x$ in the product.
\end{proposition}

\begin{proof} This is just the Taylor series for the Schwartz kernel at the
front face.
\end{proof}

\section{Cusp metrics and holomorphic families\label{A.CM}}

Although we only make we limited use of the analytic properties of cusp
metrics we recall briefly from \cite{Melrose45} and \cite{Mazzeo-Melrose4}
the basic properties of a cusp metric. Given a cusp structure on a compact
manifold with boundary, a(n exact) cusp metric is a metric on the interior
of $X$ of the form 
\begin{equation*}
g=\frac{dx^2}{x^4}+\frac h{x^2}
\end{equation*}
where $x$ is an admissible boundary defining function and $h$ is a smooth
symmetric $2$-cotensor which restricts to the boundary to be a
metric. Clearly such a metric always exists and the Laplacian, $\Lap,$ is a
cusp differential operator. Furthermore the operator $\Lap+1$ is invertible
on $\dCI(X)$ and the family $(\Lap+1)^{z/2}$ is an entire family of
c-pseudodifferential operators of complex order $z.$

\begin{lemma}\label{hop1.355} If $S\subset X$ is a closed submanifold of a
compact manifold and $F(\tau)\in I(X,S)$ is an analytic
family of $1$-step polyhomogeneous conomal functions of complex order
$-\tau$ for $\tau\in\Omega $ where $\Omega \subset\CC$ then for any
$\tau'\in\Omega$ with 
$\Re\tau'>\dim X-\dim S$ and any smooth density $\nu$ on $X$ the pairing
\begin{equation}
\langle F(\tau),\nu \rangle =\int_XF(\tau)\nu 
\label{hop1.354}\end{equation}
is defined and extends to be meromorphic in the connected component of
$\Omega $ containing $\tau'$ with only simple poles at of
$\tau=\frac14\dim X+\frac12\dim S-j,$ for $j=0,1,2,\dots.$
\end{lemma}

\noindent The position of the poles is determined by the convention for
orders of conormal distributions introduced by H\"ormander in
\cite{Hormander1}. 

In fact it is by no means necessary that the measure here be smooth,
provided it is smooth in the directions of a normal fibration to $S.$ In
particular, in the case of a manifold with corners the same result holds if
$\nu$ is smooth up to the boundary provided that $S$ is transversal to the
boundary. More generally we need to consider the case of a density which is
itself entire in another parameter in the sense that 
\begin{equation*}
\nu =\rho ^{z}\nu _0.
\end{equation*}
Here, $\nu _0$ is a fixed smooth density and $\rho$ is a defining function
for some boundary face of $X.$

\begin{lemma}\label{rbm.242} If $X$ is a compact manifold with corners and
$S\subset X$ is a closed submanifold which meets the boundary only in a
boundary hypersurface, with defining function $\rho,$ to which it is
transversal, then for $F(\tau,z)\in I(X,S)$ an analytic
family of $1$-step polyhomogeneous conomal functions of complex order
$-\tau$ for $(\tau,z)\in\Omega\times\CC$ then for any $\tau'\in\Omega$ with 
$\Re\tau'>\dim X-\dim S$ and any smooth density $\nu$ the function
\begin{equation}
\langle F(\tau),\rho^z\nu\rangle 
\label{hop1.356}\end{equation}
is defined and extends to be meromorphic in the connected component of
$\Omega\times\CC$ containing $\tau'\times\CC$ with simple poles in the two
variables at $\tau=\frac14\dim X+\frac12\dim S-j,$ for $j=0,1,2,\dots$ and
$z=-1,-2,\dots.$
\end{lemma}

\begin{proof} Choose a function $\phi\in\CIc(\RR)$ which is equal to $1$
near $0$ and has sufficiently small support. Then the function $\phi(\rho)$
localizes near the boundary hypersurfaces which $S$ meets. In the
decomposition of the pairing 
\begin{equation*}
\langle F(\tau),\rho^z\nu_0\rangle =\langle F(\tau),\phi(\rho)\rho^z\nu_0\rangle
+\langle F(\tau),(1-\phi(\rho))\rho^z\nu_0\rangle 
\end{equation*}
the second term is entire in $z$ so Lemma~\ref{rbm.242} applies with the
minor addition of an entire parameter. Taking $\phi$ to have small enough
support, allows $X$ to be replaced by a product decomposition
$H\times[0,1)$ where $H$ is the boundary hypersurface defined by $\rho$
and, in view of the transversality, $S$ can also be taken to be a product
in this set, $S=S'\times[0,1).$ Then the distribution $F(\tau)$ is a smooth
function of the coordinate $x=\rho$ with values in $I(H,S').$ The result
then follows by combining Lemma~\ref{rbm.242} in the tangential variables
with the meromorphy of the distribution $x^{z}.$
\end{proof}

\section{Cusp and b-calculi\label{A.CAbC}}

We briefly remark on the relationship between the cusp calculus considered
here and the b-calculus considered of which the K-theory is discussed in
\cite{Melrose-Nistor1}. This can be understood easily in terms of the two Lie
algebras of vector fields of which they are the respective
microlocalizations. The b-pseudodifferential operators are, in the usual
informal sense, symbolic functions of the smooth vector fields tangent to
the boundary. This is a completely natural Lie algebra of vector fields
$ \Vb(M)$ on any manifold with boundary, $M,$ with a local basis at a
boundary point, as a \ci\ module,
$$
t\frac{\pa}{\pa t},\ \frac{\pa}{\pa y_j}
$$
where $t$ is a local boundary defining function and the $y_j$ are
tangential coordinates. This should be compared with Lemma~\ref{rbm.128},
in which $x$ must define the cusp structure. Consider the smooth
transcendental map 
\begin{equation}
[0,\infty)\ni x\longmapsto t=e^{-\frac1x}\in[0,\infty).
\label{hop1.350}\end{equation}
The vector field $t\pa/\pa_t$ pulls back under this map to the smooth vector
field $x^2\pa/\pa_x.$ Thus it follows that differential operators which are
polynomials, with smooth coefficients in the vector fields $t\pa/\pa_t$
and $\pa/\pa y_j$ on $[0,1)_t\times\RR^{n-1}$ pull back to be polynomials
in the vector fields $x^2\pa/\pa_x$ and $\pa/\pa y_j.$

The map in \eqref{hop1.350} is not a diffeomorphism, so it can be
considered instead as the introduction of a singular coordinate,
$x=\frac1{\log\frac1t},$ in place of $t.$ That is, as a `transcendental blow up'
of the boundary $t=0.$ This point of view is discussed in
\cite{Hassell-Mazzeo-Melrose1} where it is shown that the new manifold with
boundary (diffeomorphic in fact to the old one) has a natural \ci\
structure. Observe that a change of $t$ in \eqref{hop1.350} to another
defining function $t'=ta(t,y),$ depending on the tangential variables,
changes $x$ to
\begin{equation*}
x'=\frac1{\log\frac1{t'}}=\frac1{\log\frac1t-\log a}=
\frac x{1-x\log a(e^{-1/x},y)}=x+O(x^2).
\end{equation*}
Thus the new manifold, with boundary defining function $x$ has a natural
cusp structure arising from the transcendental blow up. It is then
straightforward to analyze the behavior to the pseudodifferential operators.

\begin{proposition}\label{hop1.351} Under transcendental blow up of the
boundary of a compact manifold with boundary as in \eqref{hop1.350} the
b-pseudodifferential operators lift to be cusp pseudodifferential operators
and the span of the lifts over \ci\ of the blown up space is dense in the
cusp algebra.
\end{proposition}


\providecommand{\bysame}{\leavevmode\hbox to3em{\hrulefill}\thinspace}

\end{document}